\documentclass{article}
\usepackage{times}

\usepackage[square,numbers]{natbib}
\usepackage[utf8]{inputenc}
\usepackage{amsmath,amsfonts,amssymb,amsthm,mathrsfs}
\usepackage{enumitem}

\newtheorem{theorem}{Theorem}[section]
\newtheorem{lemma}[theorem]{Lemma}
\newtheorem{corollary}[theorem]{Corollary}
\newtheorem{proposition}[theorem]{Proposition}
\theoremstyle{definition}
\newtheorem{definition}[theorem]{Definition}
\newtheorem{question}[theorem]{Question}
\newtheorem{remark}[theorem]{Remark}
\date{}
\author{Danielle Witt}
\newcommand{\C}{\mathbb{C}}
\newcommand{\Z}{\mathbb{Z}}
\newcommand{\Cn}{\mathbb{C}^n}
\newcommand{\Zn}{\mathbb{Z}^n}
\newcommand{\lis}[2]{#1_1,\ldots,#1_{#2}}
\newcommand{\pows}[3]{#1_1^{#2_1}\cdots #1_{#3}^{#2_{#3}}}
\newcommand{\npows}[3]{\lVert#1_1^{#2_1}\cdots #1_{#3}^{#2_{#3}}\rVert}
\newcommand{\M}{\mathcal{M}}

\begin{document}

\title{Some results on spectra and certain norms}

\author{Danielle Witt}

\maketitle

\vspace{-0.5 in}

\renewcommand{\thefootnote}{}



\footnote{\center{ Copyright 2025, University of Illinois Urbana-Champaign}}

\renewcommand{\thefootnote}{\arabic{footnote}}
\setcounter{footnote}{0}

\begin{abstract}
Given the norms of powers $(\lVert x^n\rVert)_{n\geq 0}$ of a Banach algebra element $x$, the largest possible value of the minimum modulus on the spectrum of $x$ is determined. It is also shown that, given a Banach algebra element $x$ and a compact set $K\subset\C$ with maximum modulus no more than the spectral radius of $x$, there exists a Banach algebra element $y$ with $\lVert y^n\rVert=\lVert x^n\rVert$ for all $n\geq 0$ and spectrum equal to the union of the spectrum of $x$ and $K$. These results, along with the spectral radius formula, are generalized to the joint spectrum of several commutative Banach algebra elements. The generalization of the spectral radius formula presented gives the maximum possible joint spectrum for commutative Banach algebra elements $\lis{x}{n}$, given the norms $(\npows{x}{i}{n})_{\lis{i}{n}\geq 0}$.  
\end{abstract}

\section{Introduction}

This paper contains results regarding subsets of $\C$ that are possible spectra for a Banach algebra element $x$, given the norms $(\lVert x^n\rVert)_{n\geq 0}$, and subsets of $\Cn$ that are possible joint spectra for $n$ elements $\lis{x}{n}$ in a commutative Banach algebra, given the norms $(\npows{x}{i}{n})_{i_1,\ldots,i_n\geq 0}$. For a Banach algebra element $x$, we denote by $\sigma(x)$ the spectrum of $x$, i.e. the set of all complex numbers $\lambda$ such that $\lambda\cdot 1-x$ is not invertible. Before defining the joint spectrum, we introduce some notation. Call a linear functional $\phi$ on a complex algebra $A$ multiplicative if $\phi(xy)=\phi(x)\phi(y)$ for all $x,y\in A$.
\begin{definition}
Let $A$ be a commutative Banach algebra. Denote by $\M_A$ the set of all nonzero multiplicative linear functionals on $A$.
\end{definition}
With the Gelfand topology, the set $\M_A$ is the maximal ideal space of $A$ (see \citep[pp. 265, 268]{Rudin}).
\begin{definition}
Let $A$ be a commutative Banach algebra, and let $\lis{x}{n}\in A$. The \emph{joint spectrum} $\sigma(\lis{x}{n})$ of $\lis{x}{n}$ is given by
$$\sigma(\lis{x}{n})=\{(\phi(x_1),\ldots,\phi(x_n))\in\Cn:\phi\in\M_A\}.$$
\end{definition}

It is well known that the spectrum of an element in a Banach algebra is a nonempty compact subset of $\C$ and that the joint spectrum of $n$ elements in a commutative Banach algebra is a nonempty compact subset of $\Cn$.
When $n=1$, the joint spectrum is equal to the spectrum (see \citep[Corollary~1.2.14]{Browder}). Two necessary conditions for a sequence $(a_n)_{n\geq 0}$ of nonnegative numbers to be the sequence of norms of powers of an element in a Banach algebra are that $a_0=1$ and $a_{i+j}\leq a_ia_j$ for all $i,j\geq 0$. A fact, due to L. J. Wallen and found in \citep[pp. 50, 234]{Halmos}, is that these necessary conditions are also sufficient. 

\begin{lemma}\label{sequences}
Let $(a_n)_{n\geq 0}$ be a sequence of nonnegative  numbers. Then there exists an element $x$ in some Banach algebra with $\lVert x^n\rVert=a_n$ for all $n\geq 0$ if and only if $a_0=1$ and $a_{i+j}\leq a_ia_j$ for all $i,j\geq 0$.
\end{lemma}

An analogue to Lemma~\ref{sequences} holds in several variables and is a special case of Lemma~\ref{existence sev var}, to be proved in Section~\ref{The Joint Spectrum and Rationally Convex Connected Circled Sets, Part II}.

\begin{lemma}\label{nets}
Let $(a_{i_1\ldots i_n})_{\lis{i}{n}\geq 0}$ be an indexed collection of nonnegative
numbers. Then there exist elements $\lis{x}{n}$ in some commutative Banach algebra with $\npows{x}{i}{n}=a_{i_1\ldots i_n}$ for all $\lis{i}{n}\geq 0$ if and only if $a_{0\cdots 0}=1$ and $a_{i_1+j_1\ldots i_n+j_n}\leq a_{i_1\ldots i_n}a_{j_1\ldots j_n}$ for all
$\lis{i}{n},\lis{j}{n}\geq 0$.
\end{lemma}
In saying that a certain compact set $K\subset\C$ is a possible spectrum for a sequence of norms of powers $(a_n)_{n\geq 0}$, we shall mean that $(a_n)_{n\geq 0}$ is a sequence of nonnegative numbers such that $a_0=1$ and $a_{i+j}\leq a_ia_j$ for all $i,j\geq 0$ and that there exists an element $x$ in some Banach algebra such that $\lVert x^n\rVert=a_n$ for all $n\geq 0$ and $\sigma(x)=K$. An analogous expression will be used in several variables.

\sloppy
We denote by $r(x)$ the spectral radius of a Banach algebra element $x$, i.e. $r(x)=\sup\{|\lambda|:\lambda\in\sigma(x)\}$. Recall the spectral radius formula, that $r(x)=\lim\limits_{n\to\infty}\lVert x^n\rVert^{1/n}=\inf\limits_{n>0}\lVert x^n\rVert^{1/n}$ (see \citep[Theorem~10.13]{Rudin}). Prior work has shown that certain sequences of norms of powers give more information about the spectrum than just the spectral radius. Letting $\Gamma$ denote the unit circle, a theorem of Y. Katznelson and L. Tzafriri in \cite{Katznelson-Tzafriri} can be stated as follows. If $x$ is an element in a Banach algebra with $\sup\limits_{n>0}\lVert x^n\rVert<\infty$, then $\lVert x^n-x^{n+1}\rVert\rightarrow 0$ $(n\rightarrow\infty)$ if and only if $\sigma(x)\cap\Gamma\subset\{1\}$. (The theorem in \cite{Katznelson-Tzafriri} is stated for an operator on a Banach space, but is valid for any Banach algebra element.  This follows from the fact that every Banach algebra is isometrically isomorphic to an algebra of operators on a Banach space (see \citep[p. 230]{Rudin}).) It follows from the triangle inequality of normed vector spaces that if $x$ is a Banach algebra element, then $\lVert x^n-x^{n+1}\rVert\geq\lVert x^n\rVert-\lVert x^{n+1}\rVert$ for all $n\geq 0$. It is then not hard to use the theorem to get examples of sequences of norms of powers giving a positive spectral radius that do not permit a spectrum consisting of a single point. (For one example, let $(a_n)_{n\geq 0}$ be the sequence given by $a_{2n}=1$ and $a_{2n+1}=2$ for all $n\geq 0$.) A generalization of this theorem was proved by G. Allan and T. Ransford in \cite{Allan-Ransford}. These results give relationships between the sequence of norms of powers of a Banach algebra element and the peripheral spectrum, i.e. the intersection of the spectrum with the circle centered at the origin of radius the spectral radius.

In this paper, rather than focus on the peripheral spectrum, we give results related to the minimum modulus on the spectrum. We begin by establishing the existence of a sequence of norms of powers giving positive spectral radius that forces the spectrum to contain the origin. In Theorem~\ref{one var nonneg}, we determine the smallest annulus centered at the origin that is a possible spectrum for a given sequence of norms of powers. As a special case, we provide a characterization of the sequences of norms of powers that permit spectrum the circle centered at the origin of radius the spectral radius. It is also shown that if a certain compact set is a possible spectrum for a sequence of norms of powers $(a_n)_{n\geq 0}$, then any larger compact set contained in $\overline{D}(0,\lim\limits_{n\to\infty}a_n^{1/n})$ is also a possible spectrum for $(a_n)_{n\geq 0}$, where $\overline{D}(0,r)$ denotes the closed disc centered at the origin of radius $r\geq 0$. We provide several variable analogues of the main results. The several variable results begin with a generalization of the spectral radius formula to joint spectra, which gives the maximum possible joint spectrum for commutative Banach algebra elements $\lis{x}{n}$, given the norms $(\npows{x}{i}{n})_{\lis{i}{n}\geq 0}$. In particular, the smallest polynomially convex circled set containing the joint spectrum is determined. As Theorem~\ref{hull of spectrum}, we generalize Theorem~\ref{one var nonneg} to joint spectra, giving the intersection of all rationally convex connected circled sets that are possible joint spectra for a given collection of norms of products $(a_{\lis{i}{n}})_{\lis{i}{n}\geq 0}$.

\section{One Variable Results}\label{One Variable Results}
This section contains results on the spectrum of a single Banach algebra element. The main theorem is Theorem~\ref{one var nonneg}, which gives the largest possible value of the minimum modulus on $\sigma(x)$, given the norms of powers $(\lVert x^n\rVert)_{n\geq 0}$ of a Banach algebra element $x$. As Theorem~\ref{spectrum contains 0}, we show that there exists a sequence of norms of powers giving positive spectral radius that forces the spectrum to contain the origin. We begin with a lemma mentioned in the introduction, which shows that if a certain compact set is a possible spectrum for a sequence of norms of powers $(a_n)_{n\geq 0}$, then any larger compact set contained in $\overline{D}(0,\lim\limits_{n\to\infty}a_n^{1/n})$ is also a possible spectrum for $(a_n)_{n\geq 0}$. This is presented as Lemma~\ref{enlarging spectrum}.

Before presenting Lemma~\ref{enlarging spectrum}, we consider a few facts. If $A$ and $B$ are Banach algebras, then the direct sum of algebras $A\oplus B$ is a Banach algebra with the norm $\lVert (x,y)\rVert=\max\{\lVert x\rVert, \lVert y\rVert\}$.
The following elementary fact is not hard to verify.
\begin{lemma}\label{direct sum spectrum}
If $A$ and $B$ are Banach algebras, $x\in A$, and $y\in B$, then $\sigma((x,y))=\sigma(x)\cup\sigma(y)$, where $(x,y)$ is considered as an element of $A\oplus B$.
\end{lemma}
For a Banach algebra element $x$, we note the following alternate description of the disc $\overline{D}(0,r(x))$.
\begin{proposition}\label{disc of spectral radius} If $x$ is an element in a Banach algebra, then 
$$\overline{D}(0,r(x))=\{\lambda\in\C:|\lambda^n|\leq\lVert x^n\rVert\emph{\text{ for all }}n\geq 0\}.$$
\end{proposition}

For a compact set $K\subset\Cn$, we denote by $C(K)$ the algebra of all continuous complex-valued functions on $K$ with the supremum norm, $\lVert f\rVert=\sup\limits_{\lambda\in K}|f(\lambda)|$.

\begin{lemma}\label{enlarging spectrum}
If $x$ is an element in a Banach algebra and $K$ is a compact subset of $\overline{D}(0,r(x))$, then there exist a Banach algebra $B$ and an element $y$ of $B$ such that $\lVert y^n\rVert=\lVert x^n\rVert$ for all $n\geq 0$ and $\sigma(y)=\sigma(x)\cup K$. 
\end{lemma}
 
\begin{proof}
We assume that $K$ is nonempty, since otherwise there is nothing to prove. Let $A$ be the Banach algebra containing $x$. Let $B=A\oplus C(K)$, and let $y$ be the element $(x,z)$ of $B$, where $z$ is the inclusion of $K$ into $\C$. It follows from Proposition~\ref{disc of spectral radius} that $\lVert z^n\rVert\leq\lVert x^n\rVert$ for all $n\geq 0$. 
Therefore $\lVert y^n\rVert=\lVert x^n\rVert$ for all $n\geq 0$. It is easy to see that $\sigma(z)=K$. Lemma~\ref{enlarging spectrum} then follows from Lemma~\ref{direct sum spectrum}. 
\end{proof}

We now turn to Theorem~\ref{spectrum contains 0}.

\begin{theorem}\label{spectrum contains 0}
There exists a sequence $(a_n)_{n\geq 0}$ of positive numbers with $a_0=1$, $a_{i+j}\leq a_ia_j$ for all $i,j\geq 0$, and $\lim\limits_{n\to\infty}a_n^{1/n}=1$ such that $\sigma(x)$ must contain the origin whenever $x$ is an element in a Banach algebra with $\lVert x^n\rVert=a_n$ for all $n\geq 0$.
\end{theorem}
It is immediate from the definition of the spectrum that the spectrum of a Banach algebra element contains the origin if and only if the element is noninvertible. When $x$ is an invertible Banach algebra element, the supremum $\sup\limits_{n\geq 1}\frac{\lVert x^{n-1}\rVert}{\lVert x^n\rVert}$ is no more than $\lVert x^{-1}\rVert$. Thus, if a sequence of positive norms of powers $(a_n)_{n\geq 0}$ has $\sup\limits_{n\geq 1}\frac{a_{n-1}}{a_n}=\infty$, then the spectrum of a Banach algebra element $x$ with $\lVert x^n\rVert=a_n$ for all $n\geq 0$ must contain the origin. Theorem~\ref{spectrum contains 0} follows from the next lemma, which establishes existence of such a sequence of norms of powers giving spectral radius $1$.

\begin{lemma}\label{C-1 inf}
There exists a sequence $(a_n)_{n\geq 0}$ of positive numbers with
$a_0=1$, $a_{i+j}\leq a_ia_j$ for all $i,j\geq 0$, and $\lim\limits_{n\rightarrow\infty}a_n^{1/n}=1$ such that $\sup\limits_{n\geq 1}\frac{a_{n-1}}{a_n}=\infty$.
\end{lemma}

\begin{proof}
We define $(a_n)_{n\geq 0}$ recursively, as follows. Let $a_0=1$, $a_1=2$, and $a_2=4$. For $n\geq 3$, given $\lis{a}{n-1}$, define 
$$a_{n}=
\begin{cases}\min\{a_ia_j:i,j\geq 1\text{ and }i+j=n\}\text{ , if }a_{n-1}\leq\max\{\lis{a}{n-2}\} \\
2\text{ , otherwise}.
\end{cases}$$
We first note that
\begin{equation}\label{powers of 2}
a_n \text{ is a positive power of }2\text{ for all }n\geq 1.
\end{equation}
It is easy to see that $a_{i+j}\leq a_ia_j$ for all $i,j\geq 0$. Before turning to the remaining assertions of the lemma, we show that $(a_n)$ is unbounded. Assume, to get a contradiction, that $(a_n)$ is bounded. Then by (\ref{powers of 2}), the sequence takes on only finitely many values. We show that 
\begin{equation}\label{finitely many}
\text{for each }k\geq 1\text{, we have }a_n= 2^k\text{ for only finitely many }n,
\end{equation}
which will give a contradiction. Since $(a_n)$ takes on only finitely many values, we have $a_{n-1}>\max\{a_1,\ldots,a_{n-2}\}$ for only finitely many $n$. Consequently, only finitely many terms of $(a_n)$ are equal to $2$. Given $k\geq 1$, if $a_n\leq 2^k$ for only finitely many $n$, then $a_n\leq 2^{k+1}$ for only finitely many $n$, for if for some $N$ we have $a_n>2^k$ for all $n>N$, then $a_n>2^{k+1}$ for all $n>2N$. Statement (\ref{finitely many}) follows, and we have arrived at a contradiction. Therefore $(a_n)$ is unbounded.

Since $(a_n)$ is unbounded, we have $a_n=2$ for infinitely many $n$. Therefore $\lim\limits_{n\to\infty} a_n^{1/n}=1$. It remains to see that $\sup\limits_{n\geq 1}\frac{a_{n-1}}{a_n}=\infty$. Noting that $a_{n+1}\leq 2a_n$ for all $n\geq 0$, every nonnegative power of $2$ is achieved in the sequence. Given $j\geq 1$, if $a_k$ is the first instance of $2^{j+1}$, then $a_{k+1}=2$, and we have  
\[
\pushQED{\qed} 
\sup\limits_{n\geq 1}\frac{a_{n-1}}{a_n}\geq\frac{a_k}{a_{k+1}}=2^j.\qedhere
\popQED
\]    
\renewcommand{\qedsymbol}{}
\vspace{-\baselineskip} 
\end{proof}

\begin{remark}
The proof of Lemma~\ref{C-1 inf} remains valid if powers of $2$ are replaced by powers of $R$, for any $R>1$. 
\end{remark}

The rest of this section is dedicated to proving Theorem~\ref{one var nonneg}. Several relevant results are presented first. Theorem~\ref{one var nonneg} concerns the smallest annulus centered at the origin containing the spectrum. This annulus is necessarily closed due to the compactness of the spectrum. We use an extended definition of an annulus and regard a closed disc as a closed annulus with inner radius $0$. We first state the form of the smallest annulus centered at the origin containing the spectrum when the norms of all powers, including negative powers, if present, are known. When $x$ is an invertible Banach algebra element, this annulus has inner radius $\frac{1}{r(x^{-1})}$ by the spectral mapping theorem (see \citep[Theorem~10.28]{Rudin}).

\begin{lemma}\label{one var all powers} If $x$ is an element in a Banach algebra, then the smallest annulus centered at the origin containing $\sigma(x)$ is 
$$\begin{cases}\{\lambda\in\C:|\lambda|\leq\inf\limits_{k>0}\lVert x^k\rVert^{1/k}\}\emph{\text{ , if }}x\emph{\text{ is noninvertible}} \\
\{\lambda\in\C:\frac{1}{\inf\limits_{k>0}\lVert x^{-k}\rVert^{1/k}}\leq |\lambda|\leq\inf\limits_{k>0}\lVert x^k\rVert^{1/k}\}\emph{\text{ , if }}x\emph{\text{ is invertible}}.
\end{cases}$$
\end{lemma}

A theorem in \cite{Shields} asserts that if $T$ is a bounded unilateral weighted shift or an invertible bounded bilateral weighted shift, then $T$ has spectrum equal to the annulus from Lemma~\ref{one var all powers}. We use this to prove the next lemma. Let $\mathbb{N}_0$ denote the set of all nonnegative integers. By a unilateral weighted shift, we mean an operator $T$ on $\ell^2(\mathbb{N}_0)$ of the form $T(e_n)=\omega_ne_{n+1}$ for $n\geq 0$, where $(e_n)_{n\geq 0}$ is the standard basis for $\ell^2(\mathbb{N}_0)$ and $(\omega_n)_{n\geq 0}$ is a sequence of complex numbers. A unilateral weighted shift is necessarily noninvertible. By a bilateral weighted shift, we mean an operator $T$ on $\ell^2(\Z)$ of the form $T(e_n)=\omega_ne_{n+1}$ for $n\in\Z$, where $(e_n)_{n\in\Z}$ is the standard basis for $\ell^2(\Z)$ and $(\omega_n)_{n\in\Z}$ is a sequence of complex numbers. Part \emph{\ref{existence nonneg}} of the next lemma is Wallen's result (Lemma~\ref{sequences}) with an addendum. 
\begin{lemma}\label{existence one var}
\begin{enumerate}[label=(\roman*)]
\item\label{existence nonneg} Let $(a_n)_{n\geq 0}$ be a sequence of nonnegative numbers such that $a_0=1$ and   $a_{i+j}\leq a_ia_j$ for all $i,j\geq 0$. Then there exists an element $x$ in some Banach algebra such that $\lVert x^n\rVert=a_n$ for all $n\geq 0$ and $\sigma(x)=\{\lambda\in\C:|\lambda|\leq\inf\limits_{k>0}a_k^{1/k}\}$. 

\item\label{existence all} Let $(a_n)_{n\in\Z}$ be a sequence of positive numbers such that $a_0=1$ and $a_{i+j}\leq a_ia_j$ for all $i,j\in\Z$. Then there exists an invertible element $x$ in some Banach algebra such that
$\lVert x^n\rVert=a_n$ for all $n\in\Z$ and $\sigma(x)=\{\lambda\in\C:\frac{1}{\inf\limits_{k>0}a_{-k}^{1/k}}\leq |\lambda|\leq \inf\limits_{k>0}a_k^{1/k}\}$.

\end{enumerate}
 
\end{lemma}
Implicit in \emph{\ref{existence all}} is that $\inf\limits_{k>0}a_{-k}^{1/k}\neq 0$. The proof of Lemma~\ref{sequences} given in \citep[pp. 50, 234]{Halmos} uses a bounded unilateral weighted shift, and thus provides a proof of \emph{\ref{existence nonneg}}. The proof we provide of \emph{\ref{existence all}} is a natural modification of this proof using a bilateral weighted shift.

\begin{proof}
To see \emph{\ref{existence all}}, suppose $(a_n)_{n\in\Z}$ is a sequence of positive numbers such that $a_0=1$ and $a_{i+j}\leq a_ia_j$ for all $i,j\in\Z$. Let $T$ be the bilateral weighted shift (weighted forward shift on $\ell^2(\Z)$) given by 
$$T(e_n)=\frac{a_{n+1}}{a_n}e_{n+1}\text{ for all }n\in\Z.$$
The inequality $\frac{a_{n+1}}{a_n}\leq a_1$ holds for all $n\in\Z$, and it follows that $T$ is bounded, with $\lVert T\rVert\leq a_1$. (In fact, $\lVert T\rVert=a_1$ since $T(e_0)=a_1e_1$.)
The operator $T$ has inverse the weighted backward shift on $\ell^2(\Z)$ given by  
$$T^{-1}(e_n)=\frac{a_{n-1}}{a_n}e_{n-1}\text{ for all }n\in\Z.$$
We have
$$T^k(e_n)=\frac{a_{n+k}}{a_n}e_{n+k}\text{ for all }k,n\in\Z,$$
and it follows that $\lVert T^n\rVert= a_n$ for all $n\in\Z$. By the theorem in \cite{Shields} mentioned in the paragraph preceding the statement of the current lemma, we have $\sigma(T)=\{\lambda\in\C:\frac{1}{\inf\limits_{k>0}a_{-k}^{1/k}}\leq |\lambda|\leq \inf\limits_{k>0}a_k^{1/k}\}$. Therefore, $T$ satisfies the properties in \emph{\ref{existence all}}.
\end{proof}

We now give a definition to be used in the rest of this section. Given an indexed collection of positive numbers $(a_n)_{n\geq 0}$ such that $a_0=1$ and $a_{i+j}\leq a_ia_j$ for all $i,j\geq 0$, define, for each $i\in\Z$,
$$C_i=\sup\limits_{\substack{n\geq 0 \\ n+i\geq 0}}\frac{a_{n+i}}{a_n}.$$
Each $C_i$ lies in $(0,\infty]$. Here are two elementary properties of $(C_i)_{i\in\Z}$. Proposition~\ref{Cs finite} is used in the remainder of this section without mention.

\begin{proposition}\label{Cs equal as}
The equality $C_i=a_i$ holds for all $i\geq 0$.
\end{proposition}
\begin{proof}
Let $i\geq 0$. Since the inequality $a_{n+i}\leq a_na_i$ holds for all $n\geq 0$, we have $\frac{a_{n+i}}{a_n}\leq a_i$ for all $n\geq 0$. It follows that $C_i\leq a_i$. But also, we have $a_i=\frac{a_{0+i}}{a_0}\leq C_i$. 
\end{proof}
\begin{proposition}\label{Cs finite}
If $C_{-1}$ is finite, then $C_i$ is finite for all $i\in \Z$.
\end{proposition}
\begin{proof}
If $i<0$, then for all $n\geq 0$ such that $n+i\geq 0$, we have 
\[
\pushQED{\qed} 
\frac{a_{n+i}}{a_n}=\frac{a_{n+i}}{a_{n+i+1}}\cdots\frac{a_{n-1}}{a_n}\leq C_{-1}^{-i}.\qedhere
\popQED
\]    
\renewcommand{\qedsymbol}{}
\vspace{-\baselineskip}
\end{proof}

The next two lemmas are used to prove Theorem~\ref{one var nonneg}. The first is a simple consequence of the definition of $(C_i)_{i\in\Z}$.
\begin{lemma}\label{C lower bound norm one var}
Let $(a_n)_{n\geq 0}$ be a sequence of positive numbers such that $a_0=1$ and $a_{i+j}\leq a_ia_j$ for all $i,j\geq 0$. Suppose $x$ is an invertible Banach algebra element with $\lVert x^n\rVert=a_n$ for all $n\geq 0$. Then $C_k\leq\lVert x^k\rVert\text{ for all }k\in\Z.$
\end{lemma}
\begin{lemma}\label{Cs ineq}
Let $(a_n)_{n\geq 0}$ be a sequence of positive numbers such that $a_0=1$ and $a_{i+j}\leq a_ia_j$ for all $i,j\geq 0$. Then $C_0=1$, and if $C_{-1}$ is finite, then $C_{i+j}\leq C_iC_j\text{ for all }i,j\in\Z.$
\end{lemma}
\begin{proof}
That $C_0=1$ follows from Proposition~\ref{Cs equal as}. Suppose $C_{-1}$ is finite. 
Let $i,j\in\Z$, and let $n\geq 0$ such that $n+i+j\geq 0$. If $i\geq 0$, then $n+i\geq 0$, and we have
$$\frac{a_{n+i+j}}{a_n}=\frac{a_{n+i+j}}{a_{n+i}}\cdot\frac{a_{n+i}}{a_n}\leq C_jC_i.$$ 
If $i<0$, then $n+j>0$, and we have 
\[
\pushQED{\qed} 
\frac{a_{n+i+j}}{a_n}=\frac{a_{n+i+j}}{a_{n+j}}\cdot\frac{a_{n+j}}{a_n}\leq C_iC_j.\qedhere
\popQED
\]    
\renewcommand{\qedsymbol}{}
\vspace{-\baselineskip}
\end{proof}
 
We now present the main result of this section.

\begin{theorem}\label{one var nonneg}
Let $(a_n)_{n\geq 0}$ be a sequence of nonnegative numbers such that $a_0=1$ and
$a_{i+j}\leq a_ia_j$ for all $i,j\geq 0$.
\begin{enumerate}[label=(\roman*)]
\item\label{zeros} Suppose $a_k=0$ for some $k> 0$. If $x$ is an element in a Banach algebra with $\lVert x^n\rVert=a_n$ for all $n\geq 0$, then $\sigma(x)=\{0\}$. Furthermore, there exists an element $y$ in some Banach algebra such that $\lVert y^n\rVert=a_n$ for all $n\geq 0$.

\item\label{no zeros} Suppose each $a_n$ is positive.
Let $K$ be the annulus given by 
$$K=\begin{cases}\{\lambda\in\C:|\lambda|\leq\inf\limits_{k>0}C_k^{1/k}\}\emph{\text{ , if }}C_{-1}=\infty \\
\{\lambda\in\C:\frac{1}{\inf\limits_{k>0}C_{-k}^{1/k}}\leq |\lambda|\leq\inf\limits_{k>0}C_k^{1/k}\}\emph{\text{, if }}C_{-1}\emph{\text{ is finite}}.
\end{cases}
$$
If $x$ is an element in a Banach algebra with $\lVert x^n\rVert=a_n$ for all $n\geq 0$, then the smallest annulus centered at the origin containing $\sigma(x)$ contains $K$. Furthermore, there exists an element $y$ in some Banach algebra such that
$\lVert y^n\rVert=a_n$ for all $n\geq 0$ and $\sigma(y)=K$. 
\end{enumerate}
\end{theorem}

\begin{proof}
To see \emph{\ref{zeros}}, suppose $a_k=0$ for some $k>0$. Then $a_n=0$ for all $n\geq k$, and thus $r(x)=0$ whenever $x$ is an element in a Banach algebra with $\lVert x^n\rVert=a_n$ for all $n\geq 0$. Existence of an element $y$ in some Banach algebra such that $\lVert y^n\rVert=a_n$ for all $n\geq 0$ is given by \emph{\ref{existence nonneg}} of Lemma~\ref{existence one var}.

To see \emph{\ref{no zeros}}, suppose each $a_n$ is positive. We first show that there exists an element $y$ in some Banach algebra with $\lVert y^n\rVert=a_n$ for all $n\geq 0$ and $\sigma(y)=K$. If $C_{-1}=\infty$, then existence of an element $y$ in some Banach algebra such that $\lVert y^n\rVert=a_n$ for all $n\geq 0$ and $\sigma(y)=K$ is given by \emph{\ref{existence nonneg}} of Lemma~\ref{existence one var} and Proposition~\ref{Cs equal as}. If $C_{-1}$ is finite, then it follows from Lemma~\ref{Cs ineq} and \emph{\ref{existence all}} of Lemma~\ref{existence one var} that there exists an invertible element $y$ in some Banach algebra such that $\lVert y^n\rVert=C_n$ for all $n\in\Z$ and $\sigma(y)=K$. For this $y$, we have $\lVert y^n\rVert=a_n$ for all $n\geq 0$ by Proposition~\ref{Cs equal as}. To see the remaining assertion of \emph{\ref{no zeros}}, suppose $x$ is an element in a Banach algebra with $\lVert x^n\rVert=a_n$ for all $n\geq 0$. If $x$ is noninvertible, then by Lemma~\ref{one var all powers} and Proposition~\ref{Cs equal as}, the smallest annulus centered at the origin containing $\sigma(x)$ is $\{\lambda\in\C:|\lambda|\leq\inf\limits_{k>0}C_k^{1/k}\}$. It is clear that this annulus contains $K$. Suppose $x$ is invertible. By Lemma~\ref{C lower bound norm one var}, we have $C_{-1}\leq\lVert x^{-1}\rVert$, and thus $K=\{\lambda\in\C:\frac{1}{\inf\limits_{k>0}C_{-k}^{1/k}}\leq |\lambda|\leq\inf\limits_{k>0}C_k^{1/k}\}$. By Proposition~\ref{Cs equal as}, we have $\lVert x^k\rVert=C_k$ for all $k>0$, and it follows from Lemma~\ref{C lower bound norm one var} that $\frac{1}{\inf\limits_{k>0}\lVert x^{-k}\rVert^{1/k}}\leq\frac{1}{\inf\limits_{k>0}C_{-k}^{1/k}}$. It then follows from Lemma~\ref{one var all powers} that the smallest annulus centered at the origin containing $\sigma(x)$ contains $K$.
\end{proof}

\begin{corollary}\label{spectrum circle}
Let $(a_n)_{n\geq 0}$ be a sequence of positive numbers such that $a_0=1$, $a_{i+j}\leq a_ia_j$ for all $i,j\geq 0$, and $\lim\limits_{n\to\infty}a_n^{1/n}=r>0$. Then there exists an element $x$ in some Banach algebra with $\lVert x^n\rVert=a_n$ for all $n\geq 0$ and spectrum equal to the circle $\{\lambda\in\C:|\lambda|=r\}$ if and only if $C_{-1}$ is finite and $\inf\limits_{k>0}C_{-k}^{1/k}=\frac{1}{r}$.
\end{corollary}

Lemma~\ref{C-1 inf} gives existence of a sequence of norms of powers giving spectral radius $1$ such that $C_{-1}=\infty$. The next lemma gives, for each $0<r\leq 1$, existence of a sequence of norms of powers giving spectral radius $1$ such that $C_{-1}$ is finite and $\frac{1}{\inf\limits_{k>0}C_{-k}^{1/k}}=r$. In conjunction with Lemma~\ref{C-1 inf} (and noting Proposition~\ref{Cs equal as}), this shows that for each $0\leq r\leq 1$, the annulus $\{\lambda\in\C: r\leq |\lambda|\leq 1\}$ is a possible annulus $K$ for \emph{\ref{no zeros}} of Theorem~\ref{one var nonneg}. It follows that every closed annulus centered at the origin is a possible annulus $K$ for \emph{\ref{no zeros}} of Theorem~\ref{one var nonneg}. Indeed, if $r\geq 0$ and $R>0$, and if $K=\{\lambda\in\C: r\leq |\lambda|\leq 1\}$ for a sequence $(a_n)_{n\geq 0}$, then we have $K=\{\lambda\in\C: Rr\leq |\lambda|\leq R\}$ for the sequence $(R^na_n)_{n\geq 0}$.

\begin{lemma}
For each $0<r\leq 1$, there exists a sequence $(a_n)_{n\geq 0}$ of positive numbers with $a_0=1$, $a_{i+j}\leq a_ia_j$ for all $i,j\geq 0$, and 
$\lim\limits_{n\rightarrow\infty}a_n^{1/n}=1$ such that $C_{-1}$ is finite and $\frac{1}{\inf\limits_{k>0}C_{-k}^{1/k}}=r$. 
\end{lemma}
\begin{proof}
For $r=1$, the constant sequence given by $a_n=1$ for all $n\geq 0$ satisfies the conclusion of the lemma, for this sequence has $C_k=1$ for all $k\in\Z$.
Suppose $0<r<1$, and let $R=\frac{1}{r}$. We define $(a_n)_{n\geq 0}$ recursively, as follows. Let $a_0=1$, $a_1=R$, and $a_2=R^2$. For $n\geq 3$, given $\lis{a}{n-1}$, define 
$$a_n=\begin{cases}
R^{m-k} \text{ , if }a_{n-k}=R^m>\max\{a_1,\ldots,a_{n-k-1}\}\text{ for some }m>1\text{ and }\\
\hspace{3 in} 1\leq k\leq m-1 \\
\min\{a_ia_j:i,j\geq 1\text{ and }i+j=n\}\text{ , otherwise.}
\end{cases}$$
We note that $a_n$ is a positive power of $R$ for all $n\geq 1$. A simple modification of the proof that the sequence in Lemma~\ref{C-1 inf} is unbounded shows that $(a_n)$ is unbounded. We have $a_{n+1}\leq Ra_n$ for all $n\geq 0$, and thus, every power of $R$ is achieved in the sequence. For each $m>1$, let $n_m$ be the index of the first instance of $R^m$. We can then write, for $n\geq 3$, 
$$a_n=\begin{cases}
R^{m-k} \text{ , if }n=n_m+k\text{ for some }m>1\text{ and }1\leq k\leq m-1 \\
\min\{a_ia_j:i,j\geq 1\text{ and }i+j=n\}\text{ , otherwise.}
\end{cases}$$
To see that $a_{i+j}\leq a_ia_j$ for all $i,j\geq 0$ and that $C_{-1}$ is finite, we prove by induction on $n\geq 1$ that 
\begin{equation}\label{ai+j}
a_n\leq a_ia_j\text{ for all }n\geq 1\text{ and }i,j\geq 0\text{ such that }i+j=n\text{, and }
\end{equation} 
\begin{equation}\label{an-1}
\frac{a_{n-1}}{a_n}\leq R \text{ for all }n\geq 1.
\end{equation}
It is clear that (\ref{ai+j}) and (\ref{an-1}) are true when $n=1$ or $n=2$. Suppose $N \geq 3$ and that (\ref{ai+j}) and (\ref{an-1}) hold when $1\leq n\leq N-1$. We show that (\ref{ai+j}) and (\ref{an-1}) hold when $n=N$. If $N=n_m+k$ for some $m>1$ and $1\leq k\leq m-1$, then by definition of $(a_n)$, we have $\frac{a_{N-1}}{a_N}=R$. In this case, we also have, by the inductive assumptions, that if $i,j\geq 1$ such that $i+j=N$, then
$$a_{N}=\frac{1}{R}a_{i+j-1}\leq\frac{1}{R}a_{i-1}a_j\leq\frac{1}{R}\cdot Ra_ia_j=a_ia_j.$$
If $a_N=\min\{a_ia_j:i,j\geq 1\text{ and }i+j=N\}$, then it is clear that $a_{N}\leq a_ia_j$ for all $i,j\geq 0$ such that $i+j=N$. In this case, we have $a_N=a_ka_l$ for some $k,l\geq 1$ such that $k+l=N$. Using the inductive assumptions, we then have
$$a_{N-1}\leq a_{k-1}a_l\leq R a_ka_l=Ra_N.$$
Therefore (\ref{ai+j}) and (\ref{an-1}) are proved.

That $\lim\limits_{n\to\infty}a_n^{1/n}=1$ follows from noting that $a_{n_m+m-1}=R$ for all $m>1$. To see that $\frac{1}{\inf\limits_{k>0}C_{-k}^{1/k}}=r$, it suffices to show that $C_{-k}=R^k$ for all $k>0$. It follows from (\ref{an-1}) that $C_{-k}\leq R^k$ for all $k>0$. But also, for all $k>0$, we have
\[
\pushQED{\qed} 
C_{-k}\geq\frac{a_{n_{k+1}}}{a_{n_{k+1}+k}}=\frac{R^{k+1}}{R}=R^k.\qedhere
\popQED
\]    
\renewcommand{\qedsymbol}{}
\vspace{-\baselineskip}
\end{proof}

\section{A Generalization of the Spectral Radius Formula}\label{A Generalization of the Spectral Radius Formula}
The remaining sections concern the joint spectrum of several commutative Banach algebra elements. For the remainder of the paper, when $\lis{x}{n}$ are not explicitly defined, it is understood that $\lis{x}{n}$ are elements in an arbitrary commutative Banach algebra. In this section, we present a generalization of the spectral radius formula to joint spectra.  
In \cite{Muller}, V. M\"{u}ller determined, for $1\leq p\leq\infty$, the radius of the smallest closed $\ell^p$ ball centered at the origin containing the joint spectrum of $\lis{x}{n}$ from the norms $(\npows{x}{i}{n})_{\lis{i}{n}\geq 0}$. In this section, rather than find a radius, we determine the smallest polynomially convex circled set containing the joint spectrum. This is presented as Theorem~\ref{circled poly hull spectrum}. In Theorem~\ref{maximum joint spectrum theorem}, we show that this is the maximum possible joint spectrum for given norms of products. In Lemma~\ref{enlarging joint spectrum}, we generalize Lemma~\ref{enlarging spectrum}.

Before defining polynomial convexity and circled sets, we introduce some notation used in the remaining sections. For an element $s$ of $\Cn$, we denote by $s_k$ the $k^{th}$ coordinate of $s$, for $1\leq k\leq n$. Similarly, for an element $i$ of $\Zn$, we denote by $i_k$ the $k^{th}$ coordinate of $i$, for $1\leq k\leq n$. We denote by $S_0$ the subset of $\Zn$ of all elements with nonnegative coordinates, that is,  
$S_0=\{i\in\Zn: i_1,\ldots,i_n\geq 0\}$. For $1\leq k\leq n$, we denote by $\epsilon_k$ the element of $\Zn$ with $k^{th}$ coordinate equal to $1$ and all other coordinates equal to $0$. Addition in $\Zn$ is the usual coordinatewise addition.

\begin{definition}\label{poly convexity}
Let $K$ be a compact subset of $\Cn$. The \emph{polynomial hull} of $K$ is 
$$\widehat{K}=\{s\in\Cn:|p(s)|\leq\lVert p\rVert_K \text{ for every polynomial } p\}.$$ A compact subset $K$ of $\Cn$ is \emph{polynomially convex} if $\widehat{K}=K$.
\end{definition}
The polynomials in Definition~\ref{poly convexity} are polynomials in the complex coordinate functions $z_k:\Cn\rightarrow\C$, $1\leq k\leq n$, where $z_k$ is the projection on the $k^{th}$ coordinate. The norm $\lVert p\rVert_K$ is the supremum norm of $p$ over $K$, $\lVert p\rVert_K=\sup\limits_{s\in K}|p(s)|$.

\begin{definition}
A subset $K$ of $\Cn$ is \emph{circled} if, given  $s\in K$ and complex numbers $\lis{\alpha}{n}$ of modulus $1$, the element 
$(\alpha_1s_1,\ldots,\alpha_ns_n)$ of $\Cn$ also lies in $K$. 
\end{definition}

It follows from the maximum modulus principle that the polynomially convex circled sets in $\C$ are the closed discs centered at the origin. Before presenting Theorem~\ref{circled poly hull spectrum}, we state a relevant theorem of K. deLeeuw. In general, when an element $s$ of $\Cn$ does not lie in the polynomial hull of a compact set $K\subset\Cn$, there is a polynomial $p$ such that $|p(s)|>\lVert p\rVert_K$. The next theorem, which can be found in \cite{deLeeuw}, says that in the particular case when $K$ is circled, this polynomial can be chosen to be a monomial.

\begin{theorem}\label{Gamelin poly}
Let $K\subset\Cn$ be a compact circled set. If $s\in\Cn\setminus\widehat{K}$, then there exist nonnegative integers $\lis{i}{n}$ such that $|\pows{s}{i}{n}|>\npows{z}{i}{n}_K$.
\end{theorem}
\begin{corollary}\label{circled poly hull}
If $K$ is a compact subset of $\Cn$, then the smallest
polynomially convex circled subset of $\Cn$ containing $K$ is
$$\{s\in\Cn:|\pows{s}{i}{n}|\leq\npows{z}{i}{n}_K\emph{\text{ for all }}i\in S_0\}.$$
\end{corollary}
\begin{proof}
Let $K'=\{s\in\Cn:|\pows{s}{i}{n}|\leq\npows{z}{i}{n}_K\text{ for all }i\in S_0\}$.
It is clear that $K'$ is a polynomially convex circled set containing $K$. Suppose $K''\subset\Cn$ is a polynomially convex circled set containing $K$. If $s\in\Cn\setminus  K''$, then by Theorem~\ref{Gamelin poly}, for some $i\in S_0$, we have 
$$|\pows{s}{i}{n}|>\npows{z}{i}{n}_{K''}\geq\npows{z}{i}{n}_K,$$ so that $s\notin K'$. Therefore $K'\subset K''$. 
\end{proof}

The next lemma will be used to prove Theorem~\ref{circled poly hull spectrum}.

\begin{lemma}\label{norm equals spectral rad}
If $\lis{x}{n}$ are elements in a commutative Banach algebra, then $\npows{z}{i}{n}_{\sigma(\lis{x}{n})}=r(\pows{x}{i}{n})$ for all $i\in S_0$.
\end{lemma}
\begin{proof}
Let $A$ be the commutative Banach algebra containing $\lis{x}{n}$, and let $i\in S_0$. Then we have
\[
\pushQED{\qed}
\begin{split}
\npows{z}{i}{n}_{\sigma(\lis{x}{n})}&=\sup\{|\pows{z}{i}{n}(\phi(x_1),\ldots,\phi(x_n))|:\phi\in\M_A\} \\
&=\sup\{|\phi(x_1)^{i_1}\cdots\phi(x_n)^{i_n}|:\phi\in\M_A\}\\
&=\sup\{|\phi(\pows{x}{i}{n})|:\phi\in\M_A\} \\
&=r(\pows{x}{i}{n}).\qedhere
\end{split}
\popQED
\]    
\renewcommand{\qedsymbol}{}
\vspace{-\baselineskip}
\end{proof}

We now present the main result of this section.

\begin{theorem}\label{circled poly hull spectrum}
If $\lis{x}{n}$ are elements in a commutative Banach algebra, then the smallest polynomially convex circled subset of $\Cn$ containing $\sigma(\lis{x}{n})$ is
$$\{s\in\Cn:|\pows{s}{i}{n}|\leq\npows{x}{i}{n}\emph{\text{ for all }}i\in S_0\}.$$ 
\end{theorem}
\begin{proof}
Let $K=\{s\in\Cn:|\pows{s}{i}{n}|\leq\npows{x}{i}{n}\text{ for all }i\in S_0\}$, and let $K'=\{s\in\Cn:|\pows{s}{i}{n}|\leq r(\pows{x}{i}{n})\text{ for all }i\in S_0\}$. It follows from Corollary~\ref{circled poly hull} and Lemma~\ref{norm equals spectral rad} that the smallest polynomially convex circled subset of $\Cn$ containing $\sigma(\lis{x}{n})$ is equal to $K'$. To complete the proof, we show that $K=K'$. Since the spectral radius of any Banach algebra element is no more than its norm, we have $K'\subset K$. To see that $K\subset K'$, suppose $s\in K$. Given $i\in S_0$, we have $|s_1^{i_1k}\cdots s_n^{i_nk}|\leq\lVert x_1^{i_1k}\cdots x_n^{i_nk}\rVert$ for all $k>0$, which implies that $|\pows{s}{i}{n}|\leq\lVert (\pows{x}{i}{n})^k\rVert^{1/k}$ for all $k>0$. Thus, we have
$$|\pows{s}{i}{n}|\leq\inf\limits_{k>0}\lVert (\pows{x}{i}{n})^k\rVert^{1/k}=r(\pows{x}{i}{n})\text{ for all }i\in S_0,$$
so that $s\in K'$. Therefore $K\subset K'$.
\end{proof}

Before turning to Lemma~\ref{enlarging joint spectrum}, we present a lemma concerning the joint spectrum and the direct sum of Banach algebras, which is a generalization of Lemma \ref{direct sum spectrum}. The joint spectrum of $\lis{x}{n}$ can be written $$\sigma(\lis{x}{n})=\{s\in\Cn:1\notin(s_1\cdot 1-x_1)A+\cdots+(s_n\cdot 1-x_n)A\}.$$ 
This follows from the one-to-one correspondence between the members of $\M_A$ and the maximal ideals in $A$ that is given in \citep[p. 265]{Rudin}. 
The next lemma can be verified using this alternate formulation of the joint spectrum.

\begin{lemma}\label{joint spectrum direct sum}
Let $A$ and $B$ be commutative Banach algebras, let 
$\lis{x}{n}\in A$, and let
$\lis{y}{n}\in B$. If $(x_1,y_1),\ldots,(x_n,y_n)$ are considered as elements of $A\oplus B$, then $\sigma((x_1,y_1),\ldots,(x_n,y_n))=\sigma(\lis{x}{n})\cup\sigma(\lis{y}{n})$.
\end{lemma}

If $K$ is a nonempty compact subset of $\Cn$ and $\lis{z}{n}$ are regarded as elements of $C(K)$, then $\sigma(\lis{z}{n})=K$. This follows from the fact that the members of $\M_{C(K)}$ are exactly the evaluation functionals $\phi_s(f)=f(s)$, where $s$ ranges over $K$ (see \citep[p. 271]{Rudin}). 

\begin{lemma}\label{enlarging joint spectrum}
Let $\lis{x}{n}$ be elements in a commutative Banach algebra, and let $K$ be a compact subset of 
$\{s\in\Cn:|\pows{s}{i}{n}|\leq\npows{x}{i}{n}\emph{\text{ for all }}i\in S_0\}.$ Then there exist a commutative Banach algebra $B$ and elements $\lis{y}{n}$ in $B$ such that 
$\npows{y}{i}{n}=\npows{x}{i}{n}$ for all $i\in S_0$ and $\sigma(\lis{y}{n})=\sigma(\lis{x}{n})\cup K$. 
\end{lemma}
\begin{proof}
We assume that $K$ is nonempty, since otherwise there is nothing to prove.
Let $A$ be the commutative Banach algebra containing $\lis{x}{n}$. Let $B=A\oplus C(K)$, and for $1\leq k\leq n$, let $y_k$ be the element $(x_k,z_k)$ of $B$. It is clear that $\npows{z}{i}{n}\leq\npows{x}{i}{n}$ for all $i\in S_0$, and it follows that $\npows{y}{i}{n}=\npows{x}{i}{n}$ for all $i\in S_0$. Since $\sigma(\lis{z}{n})=K$, the lemma then follows from Lemma~\ref{joint spectrum direct sum}. 
\end{proof}

We now present Theorem~\ref{maximum joint spectrum theorem}.

\begin{theorem}\label{maximum joint spectrum theorem}
Let $(a_i)_{i\in S_0}$ be an indexed collection of nonnegative numbers such that $a_0=1$ and $a_{i+j}\leq a_ia_j$ for all $i,j\in S_0$. Then there exist elements $\lis{x}{n}$ in some commutative Banach algebra such that $\npows{x}{i}{n}=a_i$ for all $i\in S_0$ and
$\sigma(\lis{x}{n})=\{s\in\Cn:|\pows{s}{i}{n}|\leq a_i\emph{\text{ for all }}i\in S_0\}$.
\end{theorem}
\begin{proof}
Apply Lemma~\ref{nets} and apply Lemma~\ref{enlarging joint spectrum} with $K=\{s\in\Cn:|\pows{s}{i}{n}|\leq a_i\text{ for all }i\in S_0\}$.
\end{proof}

\section{The Joint Spectrum and Rationally Convex Connected Circled Sets, Part I}\label{The Joint Spectrum and Rationally Convex Connected Circled Sets, Part I}
The purpose of the remaining sections is to generalize Theorem~\ref{one var nonneg} to joint spectra. In this section, we provide a several variable generalization of Lemma~\ref{one var all powers}. In particular, as Theorem~\ref{rhull of spectrum all powers}, we determine the smallest rationally convex connected circled set containing the joint spectrum of $\lis{x}{n}$ from the norms of all defined monomials $\pows{x}{i}{n}$, where $i\in\Zn$. All compact subsets of $\C$ are rationally convex, so in $\C$, the closed annuli centered at the origin are exactly the rationally convex connected circled sets. With the exception of Lemma~\ref{rationally convex connected circled}, all results in this section are analogues of results in Section~\ref{A Generalization of the Spectral Radius Formula}.

\begin{definition}\label{rational convexity}
Let $K$ be a compact subset of $\Cn$. The \emph{rational hull} of $K$ is 
$$\widehat{K}_R=\{s\in\Cn: |f(s)|\leq\lVert f\rVert_K \text{ for every rational function }f\text{ analytic on }K\}.$$ A compact subset $K$ of $\Cn$ is \emph{rationally convex} if $\widehat{K}_R=K$.
\end{definition}

By a rational function, we mean a quotient of polynomials in the complex coordinate functions. To say that a rational function $f$ is analytic on a compact set $K\subset\Cn$ means that $f$ can be written $f=\frac{\displaystyle p}{\displaystyle q}$, where $p$ and $q$ are polynomials and $q$ has no zeroes on $K$. The norm $\lVert f\rVert_K$ is the supremum norm of $f$ over $K$, $\lVert f\rVert_K=\sup\limits_{s\in K}|f(s)|$. 

Before presenting Theorem~\ref{rhull of spectrum all powers}, we present some relevant results. Theorem~\ref{Gamelin} and Corollary~\ref{rhull of circled set} are analogues of Theorem~\ref{Gamelin poly} and Corollary~\ref{circled poly hull} for rational convexity. The next lemma is used to prove Corollary~\ref{rhull of circled set}. 

\begin{lemma}\label{rationally convex connected circled}
Let $M$ be a subset of $\mathbb{Z}^n$ containing $\epsilon_k$, for $1\leq k\leq n$, and let $(a_i)_{i\in M}$ be an indexed collection of nonnegative numbers. Then the set 
$K=\{s\in\mathbb{C}^n:|s_1^{i_1}\cdots s_n^{i_n}|\leq a_i$ \emph{for all} $i\in M\}$ is a rationally convex connected circled set. 

\end{lemma}
Requiring $M$ to contain $\epsilon_k$, for $1\leq k\leq n$, ensures that $K$ is bounded. As in Section~1, we denote the unit circle in $\C$ by $\Gamma$.

\begin{proof}
It is clear from Definition~\ref{rational convexity} that $K'\subset\widehat{K'}_R$ for any compact set $K'\subset\Cn$. To see that $K$ is rationally convex, we show that $\widehat{K}_R\subset K$. It is clear from the definition of $K$ that for each $i\in M$, the rational monomial $\pows{z}{i}{n}$ is analytic on $K$ and satisfies $\npows{z}{i}{n}_K\leq a_i$. Thus, if $s\in\widehat{K}_R$, then we have 
$$|\pows{s}{i}{n}|\leq\npows{z}{i}{n}_K\leq a_i\text{ for all }i\in M,$$
so that $s\in K$. To see that $K$ is circled, we note that if $s\in K$ and $\lis{\alpha}{n}\in\C$ with $|\alpha_1|=\cdots=|\alpha_n|=1$, then
$$|(\alpha_1s_1)^{i_1}\cdots(\alpha_ns_n)^{i_n}|=|\pows{s}{i}{n}|\leq a_i\text{ for all }i\in M,$$
so that $(\alpha_1s_1,\ldots,\alpha_ns_n)\in K$.
It remains to prove that $K$ is connected. First, define
$$B=\{(r_1,\ldots,r_n)\in\mathbb{R}^n: r_k\geq 0\text{ for each }k \text{ and }r_1^{i_1}\cdots r_n^{i_n}\leq a_i\text{ for all } i\in M\}.$$ 
The set $K$ is a continuous image of the product of $B$ and an $n$-torus, for if $F:B\times\Gamma^n\rightarrow K$ is given by $F((\lis{r}{n}),(\lis{\alpha}{n}))=(r_1\alpha_1,\ldots,r_n\alpha_n)$ for all $(\lis{r}{n})\in B$ and $(\lis{\alpha}{n})\in\Gamma^n$, then the image of $B\times\Gamma^n$ under $F$ is equal to $K$. Indeed, given $s\in K$ and writing, for each $1\leq k\leq n$, $s_k=|s_k|\alpha_k$ for some $\alpha_k\in\C$ with $|\alpha_k|=1$, we have $F((|s_1|,\ldots,|s_n|),(\alpha_1,\ldots,\alpha_n))=s$.
To prove that $K$ is connected, we show that $B$ is connected. Given $(r_1,\ldots,r_n),(s_1,\ldots,s_n)\in B$, the map $f:[0,1]\rightarrow\mathbb{R}^n$ given by
$f(t)=(r_1^{1-t}s_1^t,\ldots,r_n^{1-t}s_n^t)$ is a path in $\mathbb{R}^n$ from $(r_1,\ldots,r_n)$ to $(s_1,\ldots,s_n)$. The image of $[0,1]$ under $f$ lies in $B$, for if  $i\in M$ and $t\in [0,1]$, we have 
$$(r_1^{1-t}s_1^t)^{i_1}\cdots (r_n^{1-t}s_n^t)^{i_n}=(r_1^{i_1}\cdots r_n^{i_n})^{1-t}(s_1^{i_1}\cdots s_n^{i_n})^t\leq a_i^{1-t}a_i^t=a_i.$$ We have therefore shown that $B$ is path connected, hence connected. 
\end{proof}
The next theorem can be found in \citep[p.76]{Gamelin}.

\begin{theorem}\label{Gamelin}
Let $K\subset\Cn$ be a connected compact circled set. If 
$s\in\Cn\setminus\widehat{K}_R$, then there exist integers $i_1,\ldots,i_n$ such that the rational monomial $f(z)=z_1^{i_1}\cdots z_n^{i_n}$ is analytic on $K$ and satisfies $|f(s)|>\lVert f\rVert_K$.
\end{theorem}

\begin{corollary}\label{rhull of circled set}
If $K$ is a compact subset of $\mathbb{C}^n$, then the smallest rationally convex connected circled subset of $\Cn$ containing $K$ is 
\begin{multline*}
\{s\in\mathbb{C}^n:|s_1^{i_1}\cdots s_n^{i_n}|\leq\lVert z_1^{i_1}\cdots z_n^{i_n}\rVert_K \emph{\text{ for all }} i\in\Zn \emph{\text{ such that }} z_1^{i_1}\cdots z_n^{i_n} \\
\emph{\text{ is analytic on }}K\}.
\end{multline*} 

\end{corollary}
\begin{proof}
Let 
$K'=\{s\in\mathbb{C}^n:|s_1^{i_1}\cdots s_n^{i_n}|\leq\lVert z_1^{i_1}\cdots z_n^{i_n}\rVert_K\text{ for all } i\in\Zn \text{ such that }$\\
$z_1^{i_1}\cdots z_n^{i_n} \text{ is } \text{analytic on } K\}$.
It is clear that $K'$ contains $K$, and by Lemma~\ref{rationally convex connected circled}, $K'$ is a rationally convex connected circled set. The rest of the proof is similar to the proof of Corollary~\ref{circled poly hull}. We suppose $K''\subset\Cn$ is a rationally convex connected circled set containing $K$ and apply Theorem~\ref{Gamelin} to show that $K'\subset K''$.   
\end{proof}

The next lemma is a generalization of Lemma~\ref{norm equals spectral rad}.

\begin{lemma}\label{powers analytic lemma}
Suppose $\lis{x}{n}$ are elements in a commutative Banach algebra.
\begin{enumerate}[label=(\roman*)]
\item\label{powers analytic} For each $i\in\Zn$, the monomial $\pows{x}{i}{n}$ is defined if and only if $\pows{z}{i}{n}$ is analytic on $\sigma(\lis{x}{n})$.
\item\label{norm equals spectral rad Zn} If $i\in\Zn$ such that $\pows{x}{i}{n}$ is defined, then $\npows{z}{i}{n}_{\sigma(\lis{x}{n})}=r(\pows{x}{i}{n})$.
\end{enumerate}
\end{lemma}
\begin{proof}
The projection of $\sigma(\lis{x}{n})$ onto the $k^{th}$ coordinate is equal to $\sigma(x_k)$. Consequently, $x_k$ is invertible if and only if no element of $\sigma(\lis{x}{n})$ has $k^{th}$ coordinate equal to $0$. Statement  \emph{\ref{powers analytic}} follows. The proof of Lemma~\ref{norm equals spectral rad} is valid for all $i\in\Zn$ such that $\pows{x}{i}{n}$ is defined. This yields \emph{\ref{norm equals spectral rad Zn}}.
\end{proof}
 
We now present the main result of this section.

\begin{theorem}\label{rhull of spectrum all powers}
If $x_1,\ldots,x_n$ are elements in a commutative Banach algebra, then the smallest rationally convex 
connected circled subset of $\Cn$ containing $\sigma(\lis{x}{n})$ is 
$$\{s\in\mathbb{C}^n:|s_1^{i_1}\cdots s_n^{i_n}|\leq \lVert x_1^{i_1}\cdots x_n^{i_n}\rVert \emph{\text{ for all }} i\in\mathbb{Z}^n \emph{\text{ such that }} x_1^{i_1}\cdots x_n^{i_n} \emph{\text{ is defined}}\}.$$  
\end{theorem}
\begin{proof}
Let $K=\{s\in\mathbb{C}^n:|s_1^{i_1}\cdots s_n^{i_n}|\leq \lVert x_1^{i_1}\cdots x_n^{i_n}\rVert\text{ for all } i\in\mathbb{Z}^n\text{ such that }$\\
$x_1^{i_1}\cdots x_n^{i_n}\text{ is defined}\}$, and let $K'=\{s\in\mathbb{C}^n:|s_1^{i_1}\cdots s_n^{i_n}|\leq  r(x_1^{i_1}\cdots x_n^{i_n})\text{ for all }$\\
$i\in\mathbb{Z}^n\text{ such that } x_1^{i_1}\cdots x_n^{i_n}\text{ is defined}\}$. We conclude from Corollary~\ref{rhull of circled set} and Lemma~\ref{powers analytic lemma} that the smallest rationally convex connected circled subset of $\Cn$ containing $\sigma(\lis{x}{n})$ is equal to $K'$. To complete the proof, we show that $K=K'$. Since the spectral radius of any Banach algebra element is no more than its norm, we have $K'\subset K$. To see that $K\subset K'$, suppose $s\in K$. Given $i\in\Zn$ such that $\pows{x}{i}{n}$ is defined, we have $|s_1^{i_1k}\cdots s_n^{i_nk}|\leq\lVert x_1^{i_1k}\cdots x_n^{i_nk}\rVert$ for all $k>0$, which implies that $|\pows{s}{i}{n}|\leq\lVert (\pows{x}{i}{n})^k\rVert^{1/k}$ for all $k>0$. Thus, for all $i\in\Zn$ such that $\pows{x}{i}{n}$ is defined, we have
$$|\pows{s}{i}{n}|\leq\inf\limits_{k>0}\lVert (\pows{x}{i}{n})^k\rVert^{1/k}=r(\pows{x}{i}{n}).$$
Thus, $s\in K'$. Therefore $K\subset K'$.
\end{proof}

The next lemma is an analogue of Lemma~\ref{enlarging joint spectrum} and will be used in the next section.

\begin{lemma}\label{enlarging joint spectrum all powers}
\sloppy Let $\lis{x}{n}$ be elements in a commutative Banach algebra. Let
$M=\{i\in\Zn:\pows{x}{i}{n}\emph{\text{ is defined}}\}$, and let $K$ be a compact subset of $\{s\in\Cn:|\pows{s}{i}{n}|\leq\npows{x}{i}{n}\emph{\text{ for all }}i\in M\}$. Then there exist a commutative Banach algebra $B$ and elements $\lis{y}{n}$ in $B$ such that $\{i\in\Zn:\pows{y}{i}{n}\emph{\text{ is defined}}\}=M$, $\npows{y}{i}{n}=\npows{x}{i}{n}$ for all $i\in M$, and $\sigma(\lis{y}{n})=\sigma(\lis{x}{n})\cup K$. 
\end{lemma}
\begin{proof}
We assume that $K$ is nonempty, since otherwise there is nothing to prove. Let $A$ be the commutative Banach algebra containing $\lis{x}{n}$. Let $B=A\oplus C(K)$, and for $1\leq k\leq n$, let $y_k$ be the element $(x_k,z_k)$ of $B$. If $x_k$ is invertible, then $z_k$ is invertible, for if $x_k$ is invertible, then no element of $K$ has $k^{th}$ coordinate equal to $0$. Hence, $y_k$ is invertible if and only if $x_k$ is invertible. Therefore $\{i\in\Zn:\pows{y}{i}{n}\text{ is defined}\}=M$. It is clear that $\npows{z}{i}{n}\leq\npows{x}{i}{n}$ for all $i\in M$, and it follows that $\npows{y}{i}{n}=\npows{x}{i}{n}$ for all $i\in M$. Recalling that $\sigma(\lis{z}{n})=K$, the lemma then follows from Lemma~\ref{joint spectrum direct sum}.
\end{proof}

\section{The Joint Spectrum and Rationally Convex Connected Circled Sets, Part II}\label{The Joint Spectrum and Rationally Convex Connected Circled Sets, Part II}

In this section,  we generalize Theorem~\ref{one var nonneg} to joint spectra. In particular, we determine the largest rationally convex connected circled set that must be contained in the smallest rationally convex connected circled set containing the joint spectrum of $\lis{x}{n}$, given the norms $(\npows{x}{i}{n})_{i\in S_0}$. This is presented as Theorem~\ref{hull of spectrum}. We remark after the proof of this result that this set is the intersection of all rationally convex connected circled sets that are possible joint spectra for the given norms of products. As Question~\ref{K joint spectrum}, we leave as an open question whether this set is a possible joint spectrum for the given norms of products. A related result is presented in Theorem~\ref{spectrum norms lower bound thm}, in which we find the smallest possible value of $\npows{z}{j}{n}_{\sigma(\lis{x}{n})}$ when $j\in\Zn$ and $\pows{z}{j}{n}$ is analytic on $\sigma(\lis{x}{n})$, given the norms $(\npows{x}{i}{n})_{i\in S_0}$.

We present two main lemmas, Lemma~\ref{existence sev var} and Lemma~\ref{matching Cs}, before presenting Theorem~\ref{hull of spectrum}. The first of these is a several variable generalization of Wallen's result (Lemma~\ref{sequences}). The next lemma contains two facts that are trivial, but will be useful for the proof of Lemma~\ref{existence sev var}.

\begin{lemma}\label{lemma for existence sev var}
Let $S$ be a (possibly empty) subset of $\{1,\ldots,n\}$, and let 
$M=\{i\in\mathbb{Z}^n: i_k\geq 0\emph{\text{ for all }}k\in\{1,\ldots,n\}\setminus S\}$. Let $(a_i)_{i\in M}$ be an indexed collection of nonnegative numbers such that $a_{i+j}\leq a_ia_j$ for all $i,j\in M$. 
\begin{enumerate}[label=(\roman*)]
\item If $i,j\in M$ and $a_i=0$, then $a_{i+j}=0$. 

\item If $i,j\in M$ and $a_{i+j}\neq 0$, then $a_i\neq 0$.

\end{enumerate}
\end{lemma}

\begin{lemma}\label{existence sev var}
Let $S$ be a (possibly empty) subset of $\{1,\ldots,n\}$, and let 
$M=\{i\in\mathbb{Z}^n: i_k\geq 0\emph{\text{ for all }}k\in\{1,\ldots,n\}\setminus S\}$. Suppose $(a_i)_{i\in M}$ is an indexed collection of nonnegative numbers such that $a_0=1$ and $a_{i+j}\leq a_ia_j$ for all $i,j\in M$. Then there exist elements $x_1,\ldots,x_n$ in some commutative Banach algebra such that 
$\{i\in\mathbb{Z}^n: x_1^{i_1}\cdots x_n^{i_n}\emph{\text{ is defined}}\}=M$ and $\lVert x_1^{i_1}\cdots x_n^{i_n}\rVert=a_i$ for all $i\in M$.
\end{lemma}

For a Banach space $X$, we denote by $\mathcal{L}(X)$ the algebra of all bounded linear operators on $X$.

\begin{proof}
Given $i\in M$, let $j(i)$ be the element of $M$ given by
$$j(i)_k=\begin{cases} i_k\text{ , if } k\in S \\
0\text{ , if }k\in\{1,\ldots,n\}\setminus S.
\end{cases}$$
Noting that $-j(i)\in M$ for all $i\in M$, by Lemma~\ref{lemma for existence sev var} we have $a_{j(i)}\neq 0 \text{ for all }i\in M$. We also note that 
\begin{equation}\label{j(i+m)}
j(i+m)=j(i)+m \text{ whenever }i,m\in M \text{ and }m_k=0\text{ for all }k\in\{1,\ldots,n\}\setminus S.
\end{equation} 

We now define bounded operators $\lis{T}{n}$ on $\ell^2(M)$ and show that they satisfy the conclusion of the lemma when considered as elements of a certain commutative subalgebra of $\mathcal{L}(\ell^2(M))$. Let $(e_i)_{i\in M}$ be the standard basis for $\ell^2(M)$. For each $k\in\{1,\ldots,n\}\setminus S$, let $T_k$ be the weighted forward shift on $\ell^2(M)$ given by 
$$T_k(e_i)=\begin{cases} \frac{a_{i+\epsilon_k}}{a_i}e_{i+\epsilon_k}\text{ , if }a_i\neq 0 \\
0 \text{ , if }a_i=0
\end{cases}\text{ for all }i\in M.$$
For each $k\in S$, let $T_k$ be the weighted forward shift on $\ell^2(M)$ given by 
$$T_k(e_i)=\begin{cases} \frac{a_{i+\epsilon_k}}{a_i}e_{i+\epsilon_k}\text{ , if }a_i\neq 0 \\
\frac{a_{j(i)+\epsilon_k}}{a_{j(i)}}e_{i+\epsilon_k} \text{ , if }a_i=0
\end{cases}\text{ for all }i\in M.$$
Each $T_k$ is bounded, with $\lVert T_k\rVert= a_{\epsilon_k}$. We show that 
\begin{equation}\label{Tk commute}
\lis{T}{n}\text{ commute,}
\end{equation}
\begin{equation}\label{Tk invert}
T_k \text{ is invertible if and only if }k\in S\text{, and}
\end{equation}
\begin{equation}\label{Tk norms}
\npows{T}{i}{n}=a_i\text{ for all }i\in M.
\end{equation}

We first prove (\ref{Tk invert}). If $k\in\{1,\ldots,n\}\setminus S$, then the image of $T_k$ is contained in $\text{span}\{e_i:i\in M\text{ and }i_k>0\}$. Thus, $T_k$ is not surjective in this case. If $k\in S$, then, with some applications of Lemma~\ref{lemma for existence sev var} and (\ref{j(i+m)}), it can be seen that $T_k$ has inverse the weighted backward shift on $\ell^2(M)$ given by 
$$T_k^{-1}(e_i)=\begin{cases}\frac{a_{i-\epsilon_k}}{a_i}e_{i-\epsilon_k}\text{ , if }a_i\neq 0 \\
\frac{a_{j(i)-\epsilon_k}}{a_{j(i)}}e_{i-\epsilon_k}\text{ , if }a_i=0
\end{cases}\text{ for all }i\in M.$$

To prove (\ref{Tk commute}), we show that for all $1\leq k,l\leq n$ and $i\in M$, we have
\begin{equation}\label{TkTl}
T_kT_l(e_i)=\begin{cases}
\frac{a_{i+\epsilon_k+\epsilon_l}}{a_i}e_{i+\epsilon_k+\epsilon_l}\text{ , if }a_{i+\epsilon_k+\epsilon_l}\neq 0 \\
\frac{a_{j(i)+\epsilon_k+\epsilon_l}}{a_{j(i)}}e_{i+\epsilon_k+\epsilon_l}\text{ , if }a_{i+\epsilon_k+\epsilon_l}=0\text{ and }k,l\in S \\
0\text{ , if }a_{i+\epsilon_k+\epsilon_l}=0\text{ and either }k\notin S \text{ or }l\notin S.
\end{cases}
\end{equation}
\sloppy The first and second cases of (\ref{TkTl}) can be shown with some applications of Lemma~\ref{lemma for existence sev var}. To see the second case, also apply (\ref{j(i+m)}). To see the third case of (\ref{TkTl}), suppose $a_{i+\epsilon_k+\epsilon_l}=0$ and either $k\notin S$ or $l\notin S$. If $k\notin S$, then $T_k(e_{i+\epsilon_l})=0$, and it follows that $T_kT_l(e_i)=0$. If $k\in S$ and $l\notin S$, then an application of Lemma~\ref{lemma for existence sev var} shows that $a_{i+\epsilon_l}=0$. In this case, we then have $T_l(e_i)=0$, and therefore $T_kT_l(e_i)=0$. Thus, (\ref{TkTl}) is proved. Property (\ref{Tk commute}) follows from noting that the right side of (\ref{TkTl}) remains the same if the roles of $k$ and $l$ are reversed.

A similar proof shows that, more generally, for all $m,i\in M$, we have
\begin{equation}\label{Tpows}
\pows{T}{m}{n}(e_i)=\begin{cases}
\frac{a_{i+m}}{a_i}e_{i+m}\text{ , if }a_{i+m}\neq 0 \\
\frac{a_{j(i)+m}}{a_{j(i)}}e_{i+m}\text{ , if }a_{i+m}=0\text{ and }k\in S\text{ whenever }m_k\neq 0 \\
0\text{ , if }a_{i+m}=0\text{ and there exists }k\in\{1,\ldots,n\}\setminus S \\
\hspace{2.25 in}\text{ such that }m_k\neq 0.
\end{cases}
\end{equation}
Property (\ref{Tk norms}) follows. 

The closure of the subalgebra of $\mathcal{L}(\ell^2(M))$ generated by $\lis{T}{n}$ and the inverses $T_k^{-1}$, where $k$ ranges over $S$, is a commutative Banach algebra. The weighted shifts $\lis{T}{n}$ satisfy the conclusion of the lemma when considered as elements of this algebra.
\end{proof}

\begin{corollary}\label{rccc joint spectrum}
Let $S$ be a (possibly empty) subset of $\{1,\ldots,n\}$, and let 
$M=\{i\in\Zn: i_k\geq 0\emph{\text{ for all }}k\in\{1,\ldots,n\}\setminus S\}$. Suppose $(a_i)_{i\in M}$ is an indexed collection of nonnegative numbers such that $a_0=1$ and $a_{i+j}\leq a_ia_j$ for all $i,j\in M$. Then there exist elements $\lis{x}{n}$ in some commutative Banach algebra such that $\{i\in\Zn:\pows{x}{i}{n}\emph{\text{ is defined}}\}=M$, $\npows{x}{i}{n}=a_i$ for all $i\in M$, and $\sigma(\lis{x}{n})=\{s\in\Cn:|\pows{s}{i}{n}|\leq a_i\emph{\text{ for all }}i\in M\}$.
\end{corollary}

\begin{proof}
Apply Lemma~\ref{existence sev var} and apply Lemma~\ref{enlarging joint spectrum all powers} with
$K=\{s\in\Cn:|\pows{s}{i}{n}|\leq a_i\text{ for all }i\in M\}$, keeping in mind Theorem~\ref{rhull of spectrum all powers}. 
\end{proof}

We now give some definitions to be used in the rest of this section. Given an indexed collection $(a_i)_{i\in S_0}$ of nonnegative numbers such that $a_0=1$ and $a_{i+j}\leq a_ia_j$ for all $i,j\in S_0$, we define the following. For each $i\in\mathbb{Z}^n$, let 
$$A_i=\{m\in S_0: m+i\in S_0 \text{ and } a_m, a_{m+i}\text{ are not both }0\},$$
and define
\begin{equation*}C_i=
\begin{cases}
\sup\limits_{m\in A_i}\frac{a_{m+i}}{a_m} , \text{ if } A_i \text{ is nonempty}\\
\infty , \text{ otherwise}.
\end{cases}
\end{equation*} 
Let $S=\{k: 1\leq k\leq n \text{ and } C_{-\epsilon_k}\text{ is finite}\}$, and let 
$M=\{i\in\mathbb{Z}^n: i_k\geq 0 \text{ for all }k\in \{1,\ldots,n\}\setminus S\}$.

Each $C_i$ lies in $[0,\infty]$. It will be proved below as \emph{\ref{Cs finite sev var}} of Lemma~\ref{Cs props} that $C_i$ is finite for all $i\in M$. We note the following property of $(C_i)_{i\in\Zn}$.

\begin{proposition}\label{Cs equal as sev}
The equality $C_i=a_i$ holds for all $i\in S_0$.
\end{proposition}
\begin{proof}
Let $i\in S_0$. We have $0\in A_i$, so that $A_i$ is nonempty. If $m\in A_i$, then the inequality $a_{m+i}\leq a_ma_i$ implies that $a_m\neq 0$ and $\frac{a_{m+i}}{a_m}\leq a_i$. Therefore $C_i\leq a_i$. But also, we have $a_i=\frac{a_{0+i}}{a_0}\leq C_i$.
\end{proof}

In Theorem~\ref{hull of spectrum}, we show that the set    $K=\{s\in\mathbb{C}^n:|s_1^{i_1}\cdots s_n^{i_n}|\leq C_i$ for all $i\in M \}$ is the largest rationally convex connected circled set that must be contained in the smallest rationally convex connected circled set containing the joint spectrum of $\lis{x}{n}$ when $\npows{x}{i}{n}=a_i$ for all $i\in S_0$. When the inequality $C_{i+j}\leq C_iC_j$ holds for all $i,j\in M$, we have, from Corollary~\ref{rccc joint spectrum} and Proposition~\ref{Cs equal as sev}, existence of elements $\lis{x}{n}$ in some commutative Banach algebra such that $\npows{x}{i}{n}=a_i$ for all $i\in S_0$ and $\sigma(\lis{x}{n})=K$. In this case, Theorem~\ref{hull of spectrum} is not hard to prove. However, this inequality is not in general satisfied for all $i,j\in M$. An example can be seen with $n=2$, setting $a_{(1,0)}=2$ and $a_i=1$ for all other $i\in S_0$. We then have $M=\Z^2$ and $C_{(2,-1)}=C_{(-1,1)}=1$, while $C_{(1,0)}=2$. Lemma~\ref{matching Cs} establishes a weaker property of the collection $(C_i)_{i\in M}$ that will be used to prove Theorem~\ref{hull of spectrum}.

Before presenting Lemma~\ref{matching Cs}, we present several results used in its proof. The next lemma and its corollary give relevant properties of $(C_i)_{i\in M}$. For each integer $l$, define
$$\text{sign}(l)=\begin{cases} 1 \text{ , if } l>0 \\
0 \text{ , if }l=0 \\
-1 \text{ , if }l<0. 
\end{cases}$$

\begin{lemma}\label{Cs props}
Let $(a_i)_{i\in S_0}$ be an indexed collection of nonnegative numbers such that $a_0=1$ and $a_{i+j}\leq a_ia_j$ for all $i,j\in S_0$. Then
\begin{enumerate}[label=(\roman*)]
\item\label{Cs finite sev var} $C_i$ is finite for all $i\in M$.
\item\label{Ci+j} Suppose $i,j\in M$ satisfy
\begin{equation}\label{Ci+j crit}
m+i\in S_0\text{ or }m+j\in S_0\text{ whenever }m\in S_0\text{ and }m+i+j\in S_0.
\end{equation}
Then $C_{i+j}\leq C_iC_j$.
\end{enumerate}
\end{lemma}

\begin{proof}
We first show that
\begin{equation}\label{Ai nonempty}
A_i \text{ is nonempty for all }i\in M\text{, and }
\end{equation}
\begin{equation}\label{Ci+j Ci,j finite}
\text{ if }i,j\in M\text{ satisfy (\ref{Ci+j crit})}\text{ and }C_i\text{ and }C_j\text{ are finite, then }C_{i+j}\leq C_iC_j.
\end{equation}
To see (\ref{Ai nonempty}), we make the following observations. By definition of $S$, we have that $C_{-\epsilon_k}$ is finite for each $k\in S$. Thus, if $k\in S$, $m\in S_0$, and $a_m\neq 0$, then $a_{m+\epsilon_k}\neq 0$.
Since $a_0=1$, it follows that 
\begin{equation}\label{am neq 0 S} 
\text{if }m\in S_0\text{ and }m_k\neq 0\text{ only if }k\in S\text{, then }a_m\neq 0.
\end{equation}
Now, let $i\in M$, and let $m=(\max\{-i_1,0\},\ldots,\max\{-i_n,0\})$. Then $m$ and $m+i$ lie in $S_0$, and $m_k\neq 0$ only if $k\in S$. By (\ref{am neq 0 S}), we then have $m\in A_i$. Thus, (\ref{Ai nonempty}) is proved. 

To see (\ref{Ci+j Ci,j finite}), suppose $i,j\in M$ satisfy (\ref{Ci+j crit}) and that $C_i$ and $C_j$ are finite. The set $A_{i+j}$ is nonempty by (\ref{Ai nonempty}). We show that
$$\frac{a_{m+i+j}}{a_m}\leq C_iC_j\text{ for all }m\in A_{i+j}.$$
Let $m\in A_{i+j}$. Applying (\ref{Ci+j crit}), we assume, without loss of generality, that $m+i\in S_0$. Suppose $a_{m+i+j}\neq 0$. Then $m+i\in A_j$, and since $C_j$ and $C_i$ are finite, we have $a_{m+i}\neq 0$, $m\in A_i$, and $a_m\neq 0$. We can then write
$$\frac{a_{m+i+j}}{a_m}=\frac{a_{m+i+j}}{a_{m+i}}\cdot\frac{a_{m+i}}{a_m}\leq C_jC_i.$$
Therefore (\ref{Ci+j Ci,j finite}) is proved.

To see \emph{\ref{Cs finite sev var}}, let $i\in M$. By definition of $M$ and Proposition~\ref{Cs equal as sev}, we have that $C_{\text{sign}(i_k)\epsilon_k}$ is finite for all $1\leq k\leq n$. Also, for all $j\in M$ and $1\leq k\leq n$, we have that $j,\text{sign}(i_k)\epsilon_k$ satisfy (\ref{Ci+j crit}). Therefore, successive applications of (\ref{Ci+j Ci,j finite}) shows that
$$C_i\leq C_{\text{sign}(i_1)\epsilon_1}^{|i_1|}\cdots C_{\text{sign}(i_n)\epsilon_n}^{|i_n|}.$$
Thus, \emph{\ref{Cs finite sev var}} is proved. Statement \emph{\ref{Ci+j}} follows from \emph{\ref{Cs finite sev var}} and (\ref{Ci+j Ci,j finite}). 
\end{proof}
 
\begin{corollary}\label{Cs ineq M}
The collection $(C_i)_{i\in M}$ satisfies
$$C_{m+i}\leq C_mC_{\text{sign}(i_1)\epsilon_1}^{|i_1|}\cdots C_{\text{sign}(i_n)\epsilon_n}^{|i_n|}\text{ for all }m,i\in M.$$
\end{corollary}

We now give  some definitions to be used in the proof of Lemma~\ref{matching Cs}. Given an indexed collection $(a_i)_{i\in S_0}$ of nonnegative numbers such that $a_0=1$ and $a_{i+j}\leq a_ia_j$ for all $i,j\in S_0$, we make the following definitions. Let $p$ be the number of elements of $S$. Supposing $p>0$, and given a permutation $\lis{k}{p}$ of the elements of $S$, we define certain subsets $R_l$ of $M$ and indexed collections $(B_i^l)_{i\in M}$ of nonnegative numbers, for $0\leq l\leq p$. Define
$$R_l=\{i\in M: i_{k_{l+1}},\ldots,i_{k_{p}}\geq 0\}\text{, for }0\leq l\leq p-1\text{, and}$$
$$R_p=M.$$
We note that $R_0=S_0$ and that $R_0\subset R_1\subset\ldots\subset R_p$. Given $1\leq l\leq p$ and an indexed collection $(B_i^{l-1})_{i\in M}$ of nonnegative numbers satisfying 
\begin{enumerate}[label=(\alph*)]

\item\label{Bl-1 equal as} $B_i^{l-1}=a_i$ for all $i\in S_0$,

\item\label{Bl-1 Banach ineq} $B_{i+j}^{l-1}\leq B_i^{l-1}B_j^{l-1}$ for all $i,j\in R_{l-1}$, and 

\item\label{Bl-1 ineq M}$B_{m+i}^{l-1}\leq B_m^{l-1}(B^{l-1}_{\text{sign}(i_1)\epsilon_1})^{|i_1|}\cdots (B^{l-1}_{\text{sign}(i_n)\epsilon_n})^{|i_n|}\text{ for all }m,i\in M,$

\end{enumerate}
define 
\begin{equation}\label{Bl def}
B_i^l=\sup\limits_{m\in A_i^{l}}\frac{B_{m+i}^{l-1}}{B_m^{l-1}}\text{ for all }i\in M,
\end{equation}
where 
$$A_i^{l}=\{m\in R_{l-1}:m+i\in R_{l-1}\text{ and }B_m^{l-1},B_{m+i}^{l-1}\text{ are not both }0\}.$$
By \ref{Bl-1 equal as}, the set $A_i^l$ contains $A_i$ and is therefore nonempty by \emph{\ref{Cs finite sev var}} of Lemma~\ref{Cs props}. It follows from \ref{Bl-1 ineq M} that each $B_i^l$ is finite.  

As analogues to Proposition~\ref{Cs equal as sev}, \emph{\ref{Ci+j}} of Lemma~\ref{Cs props}, and Corollary~\ref{Cs ineq M}, we have the following properties of $(B_i^l)_{i\in M}$. These have the same proofs as the mentioned results, with the obvious substitutions.

\begin{lemma}\label{Bl analogues}
The collection $(B_i^l)_{i\in M}$ satisfies
\begin{enumerate}[label=(\roman*)]
\item\label{Bl equal Bl-1 in Rl-1}
$B_i^l=B_i^{l-1} \text{ for all } i\in R_{l-1}$,
\item\label{Bi+j}
if $i,j\in M$ satisfy
\begin{equation}\label{Bi+j crit}
m+i\in R_{l-1}\text{ or }m+j\in R_{l-1}\text{ whenever }m\in R_{l-1}\text{ and }m+i+j\in R_{l-1},
\end{equation}
then $B_{i+j}^{l}\leq B_i^lB_j^l$, and
\item\label{Bl ineq M}
$B_{m+i}^{l}\leq B_m^{l}(B^{l}_{\text{sign}(i_1)\epsilon_1})^{|i_1|}\cdots (B^{l}_{\text{sign}(i_n)\epsilon_n})^{|i_n|}\text{ for all }m,i\in M.$
\end{enumerate}
\end{lemma}

\begin{corollary}\label{Bl props}
The collection $(B_i^l)_{i\in M}$ satisfies
\begin{enumerate}[label=(\roman*)]
\item\label{Bl equal as} $B_i^{l}=a_i$ for all $i\in S_0$, and

\item\label{Bl Banach ineq} $B_{i+j}^{l}\leq B_i^{l}B_j^{l}$ for all $i,j\in R_{l}$. 
\end{enumerate}
\end{corollary}

\begin{proof}
Property \emph{\ref{Bl equal as}} follows from \ref{Bl-1 equal as} and \emph{\ref{Bl equal Bl-1 in Rl-1}} of Lemma~\ref{Bl analogues}. To see \emph{\ref{Bl Banach ineq}}, let $i,j\in R_l$. For the sake of applying \emph{\ref{Bi+j}} of Lemma~\ref{Bl analogues}, suppose $m\in R_{l-1}$ such that $m+i+j\in R_{l-1}$. It is not hard to see that if $i_{k_l}\geq 0$, then $m+i\in R_{l-1}$, and if $i_{k_l}<0$, then $m+j\in R_{l-1}$. Therefore $i,j$ satisfy (\ref{Bi+j crit}), and by \emph{\ref{Bi+j}} of Lemma~\ref{Bl analogues}, we have $B_{i+j}^l\leq B_i^lB_j^l$.
\end{proof}

We now set $B_i^0=C_i$ for all $i\in M$ and define indexed collections $(B_i^l)_{i\in M}$ of nonnegative numbers, for $0\leq l\leq p$, recursively using (\ref{Bl def}). By Proposition~\ref{Cs equal as sev} and Corollary~\ref{Cs ineq M}, the collection $(B_i^0)_{i\in M}$ satisfies \ref{Bl-1 equal as}, \ref{Bl-1 Banach ineq}, and \ref{Bl-1 ineq M}. This recursive definition is valid by Corollary~\ref{Bl props} and \emph{\ref{Bl ineq M}} of Lemma~\ref{Bl analogues}. Noting Proposition~\ref{Cs equal as sev}, we have the following remark direct from the definition of $(B_i^1)_{i\in M}$.

\begin{remark}\label{B1 equal C}
The equality $B_i^1=C_i$ holds for all $i\in M$.
\end{remark}

The next two lemmas give relevant properties of the collections $(B_i^l)_{i\in M}$, $0\leq l\leq p$. 

\begin{lemma}\label{Bl-1 leq Bl}
The inequality $B_i^{l-1}\leq B_i^l$ holds for all $1\leq l\leq p$ and $i\in M$.
\end{lemma}

\begin{proof}
When $l=1$, the lemma holds by Remark~\ref{B1 equal C}. Suppose $2\leq l\leq p$, and let $i\in M$. Applying \emph{\ref{Bl equal Bl-1 in Rl-1}} of Lemma~\ref{Bl analogues} and noting that $A_i^{l-1}\subset A_i^{l}$, we have
\[
\pushQED{\qed}
B_i^{l-1}=\sup\limits_{m\in A_i^{l-1}}\frac{B_{m+i}^{l-2}}{B_m^{l-2}}=\sup\limits_{m\in A_i^{l-1}}\frac{B_{m+i}^{l-1}}{B_m^{l-1}}\leq\sup\limits_{m\in A_i^l}\frac{B_{m+i}^{l-1}}{B_m^{l-1}}=B_i^l.\qedhere
\popQED
\]     
\renewcommand{\qedsymbol}{}
\vspace{-\baselineskip}
\end{proof}

\begin{lemma}\label{Bl equal Bl-1}
Suppose $2\leq l\leq p$. If $i\in M$ such that $i_{k_{l-1}}\leq 0$, then $B_i^l=B_i^{l-1}$.
\end{lemma}

\begin{proof}
Let $2\leq l\leq p$, and suppose $i\in M$ such that $i_{k_{l-1}}\leq 0$. By Lemma~\ref{Bl-1 leq Bl}, it suffices to show that $B_i^l\leq B_i^{l-1}$. For all $j\in R_{l-1}$, we have that $i,j$ satisfy (\ref{Bi+j crit}), with $l-1$ in place of $l$. Indeed, if $j\in R_{l-1}$, $m\in R_{l-2}$, and $m+i+j\in R_{l-2}$, it is not hard to see that since $i_{k_{l-1}}\leq 0$, we have $m+j\in R_{l-2}$. We can then apply \emph{\ref{Bi+j}} of Lemma~\ref{Bl analogues}, with $l-1$ in place of $l$, to see that $B_{j+i}^{l-1}\leq B_j^{l-1}B_i^{l-1}$ for all $j\in A_i^l$. It follows that $B_i^l\leq B_i^{l-1}$. 
\end{proof}

We now present Lemma~\ref{matching Cs}.

\begin{lemma}\label{matching Cs}
Let $(a_i)_{i\in S_0}$ be an indexed collection of nonnegative numbers such that $a_0=1$ and
$a_{i+j}\leq a_ia_j$ for all $i,j\in\ S_0$. Then for each $m'\in M$, there exists an indexed collection $(B_i)_{i\in M}$ of nonnegative numbers such that $B_{i+j}\leq B_iB_j$ for all $i,j\in M$, $B_i=a_i$ for all $i\in S_0$, $C_i\leq B_i$ for all $i\in M$, and $B_{m'}=C_{m'}$.
\end{lemma}

\begin{proof}
If $S$ is empty, then $M=S_0$, and, noting Proposition~\ref{Cs equal as sev}, the lemma is trivial. Suppose $S$ is nonempty, and let $m'\in M$. Let $N$ be the number of negative coordinates of $m'$. We note that $N\leq p$. If $N\geq 1$, we let $\lis{k}{p}$ be a permutation of the elements of $S$ such that $\lis{k}{N}$ are the distinct elements $k$ of $\{1,\ldots,n\}$ such that $m'_k<0$. If $N=0$, we let $\lis{k}{p}$ be any permutation of the elements of $S$. Define $(B_i^l)_{i\in M}$, for $0\leq l\leq p$, as in the paragraph following Corollary~\ref{Bl props}, and let
$$B_i=B_i^p\text{ for all }i\in M.$$
We show that $(B_i)_{i\in M}$ satisfies the properties of the lemma.
 
By Corollary~\ref{Bl props}, we have that $B_i=a_i$ for all $i\in S_0$ and $B_{i+j}\leq B_iB_j$ for all $i,j\in M$. It follows from Lemma~\ref{Bl-1 leq Bl} that $C_i\leq B_i$ for all $i\in M$. It remains to show that $B_{m'}=C_{m'}$. First, define
$$S_l=\{i\in M : i_{k_1}\leq 0,\ldots,i_{k_l}\leq 0,i_{k_{l+1}}\geq 0,\ldots,i_{k_p}\geq 0\}\text{, for }1\leq l\leq p-1\text{, and}$$
$$S_p=\{i\in M : i_{k_1},\ldots,i_{k_p}\leq 0\}.$$
We note that
$S_l\subset R_l\text{ for all }0\leq l\leq p$.
By the choice of $\lis{k}{p}$, we have that
\begin{equation}\label{m' in SN}
m'\in S_N.
\end{equation}
By Remark~\ref{B1 equal C} and Lemma~\ref{Bl equal Bl-1}, we have
$$B_i^q=B_i^{q-1}\text{ whenever }1\leq q\leq l\leq p\text{ and }i\in S_l.$$
It follows that 
$$B_i^l=C_i\text{ for all }0\leq l\leq p\text{ and }i\in S_l.$$
It then follows from (\ref{m' in SN}) and \emph{\ref{Bl equal Bl-1 in Rl-1}} of Lemma~\ref{Bl analogues} that $B_{m'}=C_{m'}$.
\end{proof}

\begin{remark}\label{Bs equal Cs in Sl}
The proof of Lemma~\ref{matching Cs} actually shows that $B_i=C_i$ for all $i\in\bigcup\limits_{l=0}^{p}S_l$.

\end{remark}

The next lemma will also be used to prove Theorem~\ref{hull of spectrum}.

\begin{lemma}\label{Cs signif}
Let $(a_i)_{i\in S_0}$ be an indexed collection of nonnegative numbers such that $a_0=1$ and $a_{i+j}\leq a_ia_j$ for all $i,j\in S_0$. Suppose $\lis{x}{n}$ are elements in a commutative Banach algebra with $\npows{x}{i}{n}=a_i$ for all $i\in S_0$. If $i\in\Zn$ such that $\pows{x}{i}{n}$ is defined, then
\begin{enumerate}[label=(\roman*)]
\item\label{i in M} $i\in M$, and
\item\label{C lower bound norm} $C_i\leq\npows{x}{i}{n}$.
\end{enumerate}
\end{lemma}

\begin{proof}
Let $i\in\Zn$ such that $\pows{x}{i}{n}$ is defined. To prove \emph{\ref{i in M}}, we show that $C_{-\epsilon_k}$ is finite if $i_k<0$. For each $1\leq k \leq n$, the set $A_{-\epsilon_k}$ is nonempty, since $\epsilon_k\in A_{-\epsilon_k}$. If $i_k<0$, then $x_k$ is invertible, and it is easy to see that $C_{-\epsilon_k}\leq\lVert x_k^{-1}\rVert$. Thus, \emph{\ref{i in M}} is proved. It then follows from \emph{\ref{Cs finite sev var}} of Lemma~\ref{Cs props} that $A_i$ is nonempty, and it is then easy to see \emph{\ref{C lower bound norm}}.
\end{proof}

We now present the main result of this section.

\begin{theorem}\label{hull of spectrum}
Let $(a_i)_{i\in S_0}$ be an indexed collection of nonnegative numbers such that $a_0=1$ and
$a_{i+j}\leq a_ia_j$ for all $i,j\in S_0$. Let 
$$K=\{s\in\mathbb{C}^n:|s_1^{i_1}\cdots s_n^{i_n}|\leq C_i \emph{\text{ for all }} i\in M \}.$$ Then 
\begin{enumerate}[label=(\roman*)]

\item\label{K char} $K$ is a rationally convex connected circled set, and 

\item\label{K prop} whenever $x_1,\ldots,x_n$ are elements in a commutative Banach algebra with
$\lVert x_1^{i_1}\cdots x_n^{i_n}\rVert=a_i$ for all $i\in S_0$, the smallest rationally convex connected circled subset of $\Cn$ containing $\sigma(x_1,\ldots,x_n)$ contains $K$.

\end{enumerate}
Furthermore, $K$ is the largest subset of $\mathbb{C}^n$ for which \ref{K prop} holds. 

\end{theorem}

To say that $K$ is the largest subset of $\Cn$ for which \emph{\ref{K prop}} holds means that if \emph{\ref{K prop}} holds when $K$ is replaced by a set $K'\subset\Cn$, then $K'\subset K$.

\begin{proof}
Statement \emph{\ref{K char}} follows from Lemma~\ref{rationally convex connected circled}. Statement \emph{\ref{K prop}} follows from Lemma~\ref{Cs signif} and Theorem~\ref{rhull of spectrum all powers}. 
It remains to see that $K$ is the largest subset of $\Cn$ for which \emph{\ref{K prop}} holds.
Let $\mathscr{E}$ be the collection of all sets $E\subset\mathbb{C}^n$ such that $E$ is the smallest  rationally convex connected circled subset of $\Cn$ containing $\sigma(x_1,\ldots,x_n)$ for some commutative Banach algebra elements $x_1,\ldots,x_n$ with
$\lVert x_1^{i_1}\cdots x_n^{i_n}\rVert=a_i$ for all $i\in S_0$. We show that $K$ is equal to the intersection $\bigcap\limits_{E\in\mathscr{E}}E$.
That $K\subset\bigcap\limits_{E\in\mathscr{E}}E$ follows from \emph{\ref{K prop}}. To see that $\bigcap\limits_{E\in\mathscr{E}}E\subset K$, let $\lambda\in\bigcap\limits_{E\in\mathscr{E}}E$. By Lemma~\ref{matching Cs}, Lemma~\ref{existence sev var}, and Theorem~\ref{rhull of spectrum all powers}, for each $m\in M$, there is a member of $\mathscr{E}$ of the form $\{s\in\Cn:|\pows{s}{i}{n}|\leq B_i\text{ for all }i\in M\}$ for some indexed collection $(B_i)_{i\in M}$ of nonnegative numbers such that $B_m=C_m$. It follows that $|\pows{\lambda}{i}{n}|\leq C_i$ for all $i\in M$, so that $\lambda\in K$. Thus, $\bigcap\limits_{E\in\mathscr{E}}E\subset K$.
\end{proof}

By Theorem~\ref{rhull of spectrum all powers} and Corollary~\ref{rccc joint spectrum}, each member of the collection $\mathscr{E}$ in the proof of Theorem~\ref{hull of spectrum} is the joint spectrum of some commutative Banach algebra elements $\lis{x}{n}$ with $\npows{x}{i}{n}=a_i$ for all $i\in S_0$. 
Thus, the proof of Theorem~\ref{hull of spectrum} shows that $K$ is the intersection of all rationally convex connected circled sets that are possible joint spectra for the norms of products $(a_i)_{i\in S_0}$. Since the inequality $C_{i+j}\leq C_iC_j$ is not in general satisfied for all $i,j\in M$, it is unclear if $K$ is a possible joint spectrum for $(a_i)_{i\in S_0}$. This is left as an open question. 

\begin{question}\label{K joint spectrum}
If $(a_i)_{i\in S_0}$ is an indexed collection of nonnegative numbers with $a_0=1$ and $a_{i+j}\leq a_ia_j$ for all 
$i,j\in S_0$ and $K$ is defined as in Theorem~\ref{hull of spectrum}, do there exist elements $\lis{x}{n}$ in some commutative Banach algebra such that $\npows{x}{i}{n}=a_i$ for all $i\in S_0$ and $\sigma(\lis{x}{n})=K$? 
\end{question}

We now present Theorem~\ref{spectrum norms lower bound thm}.
Recall that for $i\in\Zn$, the monomial $\pows{x}{i}{n}$ is defined if and only if $\pows{z}{i}{n}$ is analytic on $\sigma(\lis{x}{n})$. This is \emph{\ref{powers analytic}} of Lemma~\ref{powers analytic lemma}. For $k>0$ and $i\in\Zn$, we denote by $ki$ the element $(ki_1,\ldots,ki_n)$ of $\Zn$.
\begin{theorem}\label{spectrum norms lower bound thm}
Let $(a_i)_{i\in S_0}$ be an indexed collection of nonnegative numbers such that $a_0=1$ and $a_{i+j}\leq a_ia_j$ for all $i,j\in S_0$. 
\begin{enumerate}[label=(\roman*)]
\item\label{spectrum norms lower bound} \sloppy Suppose $\lis{x}{n}$ are elements in a commutative Banach algebra with $\npows{x}{i}{n}=a_i$ for all $i\in S_0$. If $i\in\Zn$ and $\pows{x}{i}{n}$ is defined, then $i\in M$ and $\npows{z}{i}{n}_{\sigma(\lis{x}{n})}\geq\inf\limits_{k>0}C_{ki}^{1/k}$.
\item\label{lower bound reached} For each $j\in M$, there exist elements $\lis{y}{n}$ in some commutative Banach algebra such that $\npows{y}{i}{n}=a_i$ for all $i\in S_0$, $\pows{y}{j}{n}$ is defined, and $\npows{z}{j}{n}_{\sigma(\lis{y}{n})}=\inf\limits_{k>0}C_{kj}^{1/k}$.
\end{enumerate}
\end{theorem}

\begin{proof}

Statement \emph{\ref{spectrum norms lower bound}} follows from Lemma~\ref{powers analytic lemma} and Lemma~\ref{Cs signif}. Statement \emph{\ref{lower bound reached}} follows from Lemma~\ref{powers analytic lemma}, Lemma~\ref{matching Cs} with Remark~\ref{Bs equal Cs in Sl}, and Lemma \ref{existence sev var}.
\end{proof}

\subsection*{Acknowledgements}
This paper consists of results from the author's dissertation. The author thanks her advisor, Alexander Izzo, for directing the research that led to this paper and for providing many helpful suggestions. The author also thanks Anthony O'Farrell for mention of the work of Katznelson and Tzafriri in \cite{Katznelson-Tzafriri}.

\vspace{0.25 in}

\emph{Email address: danwitt@bgsu.edu}


\begin{thebibliography}{10}

\bibitem{Allan-Ransford} G. Allan and T. Ransford, \emph{Power-dominated elements in a Banach algebra}, Studia Mathematica 94 (1989), 63-79.

\bibitem{Browder} A. Browder, \emph{Introduction to Function Algebras}, Benjamin, New York, 1969.

\bibitem{deLeeuw} K. deLeeuw, \emph{A type of convexity in the space of $n$ complex variables}, Transactions of the American Mathematical Society 83 (1956), 193--204.

\bibitem{Gamelin} T.W. Gamelin, \emph{Uniform Algebras}, $2^{nd}$ edition, Chelsea, New York, 1984.

\bibitem{Halmos} P.R. Halmos, \emph{A Hilbert Space Problem Book}, $2^{nd}$ edition, Springer-Verlag, New York, 1982. 

\bibitem{Katznelson-Tzafriri} Y. Katznelson and L. Tzafriri, \emph{On power bounded operators}, Journal of Functional Analysis 68 (1986), 313--328.    


\bibitem{Muller} V. M\"{u}ller, \emph{On the joint spectral radius}, Annales Polonici Mathematici 66 (1997), 173--182.

\bibitem{Rudin} W. Rudin, \emph{Functional Analysis}, McGraw-Hill, New York, 1973. 

\bibitem{Shields} A.L. Shields, \emph{Weighted shift operators and analytic function theory}, in: Topics in Operator Theory, Mathematical Surveys vol. 13, American Mathematical Society, Providence, RI, 1974, 49--128.

\end{thebibliography}
\end{document}